\documentclass[11pt]{amsart}
\usepackage{amssymb,hyperref,times}
\newtheorem{prop}{Proposition}[subsection]

\title{Floyd maps for relatively hyperbolic groups}

\author
{Victor Gerasimov\\ \\ $\text{\today}$}
\address{Victor Guerassimov, Departamento de Matem\'atica,
Universidade Federal de Minas Gerais,
Av. Ant\^onio Carlos, 6627/
CEP: 31270-901. Caixa Postal: 702
Belo Horizonte, MG,
Brazil}
\email{victor@mat.ufmg.br}

\thanks{This work was partially supported by Max Plank Institut of Mathematics, Bonn, Germany
and by Universit\'e Lille 1, France}


\begin{document}
%\large
%\LARGE
%\Large
\def\d{\boldsymbol\delta}
\def\eod{\hfill$\square$}
\let\Cal\mathcal
% [abstract] @tag_4A80EE9B:
\begin{abstract}
%\LARGE
%\Large
Let $\d_{\Cal S,\lambda}$ denote
the Floyd metric on a discrete group $G$ generated by a finite set $\Cal S$
with respect to the scaling function $f_n{=}\lambda^n$ for a positive $\lambda{<}1$.
We prove that if $G$ is relatively hyperbolic with respect to a collection $\mathcal P$ of subgroups
then there exists $\lambda$ such that
the identity map $G\to G$ extends to a continuous equivariant map
from the completion with respect to $\d_{\Cal S,\lambda}$
to the Bowditch completion of $G$ with respect to $\mathcal P$.

In order to optimize the proof and the usage of the map theorem
we propose two new definitions of relative hyperbolicity
equivalent to the other known definitions.

In our approach some ``visibility'' conditions in graphs are essential.
We introduce a class of ``visibility actions''
that contains the class of relatively hyperbolic actions.
The convergence property still holds for the visibility actions.

Let a locally compact group $G$ act on a compactum $\Lambda$ with convergence property
and on a locally compact Hausdorff space $\Omega$ properly and cocomactly.
Then the topologies on $\Lambda$ and $\Omega$ extend uniquely
to a topology on the direct union $T=\Lambda{\sqcup}\Omega$
making $T$ a compact Hausdorff space such that the action $G{\curvearrowright}T$
has convergence property.
We call $T$ the \it attractor sum \rm of $\Lambda$ and $\Omega$.

%inculde ``locally compact groups'' and ``continuous actions'' to the subject classification
\end{abstract}
\date{\today}
\maketitle
\def\colontitle{V.~Gerasimov, {\sl Floyd maps for relatively hyperbolic groups}\quad(\today)}
\markboth{\colontitle}{\colontitle}
% [myDefs] @tag_4A80F136:
%myDefs
%
\hsize400pt
\leftskip-50pt
\def\tupl#1{\left(#1\right)}
\let\on\curvearrowright
\let\b\mathbf
\let\s\mathsf
\let\g\mathfrak
\let\Cal\mathcal
\let\goth\mathfrak
\font\cFont=cmr9%
\def\note#1{\{{\cFont #1}\}}%
% [Introduction] @tag_4DEB8D8F(
\section{Introduction}
% [floyd thm] @tag_4DEB8D9F:
\subsection{Floyd Construction}
In \cite{Fl80} W. J. Floyd introduced a class of metrics on finitely generated
groups obtained by the ``conformal scaling'' of the word metric.
Namely he regards a group $G$ as the vertex set of a
locally finite \bf metric \rm Cayley graph,
where the length of an edge
depends on the word distance from the origin.
The function $f:n\mapsto f_n{=}($the new length of an edge of word distance $n$ from $1)$
is called the \it scaling function\rm.
Under certain conditions on $f$
the Cauchy completion $\s{Fl}_fG$ with respect
to the new path metric is compact and $G$ acts on $\s{Fl}_fG$ by bi-lipschitz homeomorphisms.

Floyd proved that, for any finitely generated
geometrically finite Kleinian group $G$
and for the scaling function $f:n\mapsto{1\over n^2+1}$, every orbit map
$G{\ni}g\mapsto g\goth p{\in}\mathbb H^3$ extends by continuity to the
\it Floyd map \rm$\phi_f:\s{Fl}_fG\to\overline{\mathbb H^3}=\mathbb H^3{\cup}\mathbb S^2$ that takes
the \it Floyd boundary \rm$\partial_fG{=}\s{Fl}_fG{\setminus}G$
onto the ``limit set'' $\b\Lambda G$.

Floyd also calculated the ``kernel'' of the boundary map $\phi_f|_{\partial_fG}$:
the map is one-to-one except
for the preimages of the parabolic points of rank one where it is two-to-one.

Our purpose is to generalize Floyd's result replacing
`geometrically finite Kleinian' by `relatively hyperbolic' (r.h. for short).
It this paper we prove the existence of the Floyd map.
In the next papers \cite{GP09}, \cite{GP10} we describe its kernel.
% [rel hyperbolicity] @tag_4DEB8E4C:
\subsection{Relative hyperbolicity}
Our strategy depends on the choice of the
initial definition of the relative hyperbolicity.
In \cite{Hr10} some relations between various definitions are discussed.
The following two main definitions, the ``geometric'' and the ``dynamical'', reflect two aspects of the subject.

A graph is said to be \it fine \rm\cite{Bo97} if
each set of arcs of bounded length with fixed
endpoints is finite (see \ref{grphs}).

\def\deff#1{%
%\box255=%
\dimen255\hsize\relax\advance\dimen255 by -64pt%
\vskip3pt\hbox{\strut\vrule\kern3pt\vrule\kern60pt%
\vbox{\hsize\dimen255\noindent\strut\bf Definition\rm #1\relax\strut}}\vskip3pt}

\deff{ $\s{RH_{fh}}$ \rm(`fh' stands for `fine hyperbolic').
An action of a group $G$ on a connected graph $\Gamma$ is said to be \it relatively hyperbolic \rm
 if $\Gamma$ is $\delta$-hyperbolic and fine,
the action $G{\on}\Gamma^1{=}\{$the edges of $\Gamma\}$
is proper (i.e, the stabilizers of the edges are finite), cofinite (i.e, $\Gamma^1/G$ is finite),
and non-parabolic (no vertex is fixed by the whole $G$).}
\vskip3pt
 A group $G$ is
said to be \it relatively hyperbolic
with respect to \rm a finite collection $\Cal P$ of infinite subgroups
if it possesses an r.h. action $G{\on}\Gamma$ such that $\Cal P$ is a set of representatives
of the orbits of the stabilizers of the vertices of infinite degree. 
\vskip3pt
Recall that a continuous action of a locally compact group $G$ on a locally compact
Hausdorff space $L$
is said to be \it proper \rm(see \ref{charProp}) if for every compact set $K{\subset}L{\times}L$ the
set $\{(g,\goth p){\in}G{\times}L:(\goth p,g\goth p){\in}K\}$ is compact.
Proper actions of discrete groups are also called \it properly discontinuouos\rm.
If $L$ is also discrete then `proper' means `the stabilizers of points are finite'.

An action is said to be \it proper on triples\rm, if the induces action on the
space of subsets of cardinality 3 is proper.

\deff{ $\s{RH_{32}}$ \rm(`$32$' means `3-proper and 2-cocompact').
An action of a discrete group $G$
by auto-homeomorphisms of a compactum $T$
is said to be \it relatively hyperbolic \rm if
it is proper on triples, cocompact on pairs and has at least two ``limit points'' (see \ref{limitSet}).}

It is less easy to restore the collection of ``parabolic'' subgroups starting from this
definition. Let us call a point $\goth p{\in}T$ \it parabolic\rm, if its stabilizer $H$
in $G$ is infinite and acts cocompactly on $T{\setminus}\{\goth p\}$.

A group $G$ is
said to be \it relatively hyperbolic
with respect to \rm a finite collection $\Cal P$ of infinite subgroups
if it possesses an r.h. (in the sense of $\s{RH}_{32}$) action $G{\on}T$ such that
$\Cal P$ is a set of representatives
of the stabilizers of parabolic (see \ref{cones}) points.
\vskip3pt
In both cases (as well as in other equivalent definitions) an r.h. \bf group \rm with respect to $\Cal P$ is defined by means of an r.h. \bf action \rm
with a specified set of ``parabolic'' subgroups which is nontrivial in some sense (``nonparabolic'').
The notion of r.h.action seems to be more fundamental than the notion of r.h.group.

At the present moment no simple proof of `$\s{RH_{fh}{\Leftrightarrow}RH_{32}}$' is known.
To obtain either of the implications one should interpret
an r.h. action as an action of the other type.
So, for `$\Rightarrow$', one has to
construct a compactum $T$ acted upon by $G$
with the desired properties,
and, for `$\Leftarrow$' one has to
find at least one fine hyperbolic $G$-graph (i.e. a graph endowed with an action of $G$) out of the topological information.

This paper is intended to facilitate the translation between the ``geometric'' and the ``dynamical''
languages in both directions.
In \cite{GP10} and \cite{GP11} we will advance in this program.

Our proof of the Floyd map theorem \ref{mapThrm} requires some information
derivable from \bf both \rm of the characteristic properties $\s{RH}_*$.
Several attempts to find a short self-contained proof
starting from just one of them failed.
This is a serious shortcoming since Floyd theorem
proved to be
useful in the development of
the theory from either of $\s{RH}_*$ and in its applications.
To resolve this problem we propose
two intermediate definitions, $\s{RH_{ah}}$ and $\s{RH_{pd}}$
(equivalent to the other $\s{RH}_*$)
of ``minimal distance'' between geometry and dynamics.
The most difficult implication `$\s{RH_{32}{\Rightarrow}RH_{fh}}$'
splits into problems of different nature:
$$\s{RH_{32}{\Rightarrow}RH_{pd}{\Rightarrow}RH_{ah}{\Rightarrow}RH_{fh}}.$$
The first implication is some finiteness problem,
the second is a deduction `dynamics$\Rightarrow$geometry'
and the third is an easy geometry.
% [alt-hypebolicity] @tag_4DFA203A:
\subsection{Alternative hyperbolicity}
We call a connected graph $\Gamma$
\it alternatively hyperbolic \rm if
for every edge $e$ of $\Gamma$ there exists
a finite set $F{\subset}\Gamma^1$ such that every geodesic triangle containing
$e$ on a side also contains an edge in $F$ on another side.

\deff{ $\s{RH_{ah}}$\rm.
An action of a group $G$ on a connected graph $\Gamma$ is said to be \it relatively hyperbolic \rm
if $\Gamma$ is alternatively hyperbolic and
the action $G{\on}\Gamma^1$
is proper cofinite and no vertex is fixed by $G$.}
% [divider] @tag_4DFCD6DB:
\subsection{Perspective divider}\label{iDivider}
We express relative hyperbolicity as a uniform structure on a group $G$
or on a ``connected $G$-set''.
Recall that a \it uniformity \rm
on a set $M$ is a filter $\Cal U$ on the set $M^2{=}M{\times}M$ whose elements are called \it entourages\rm.
Each entourage $\b u$ should contain the diagonal $\b\Delta^2M$ and a symmetric entourage $\b v$ such that
$\b v^2{\subset}\b u$. We often regard an entourage as a set of non-ordered pairs.

Let $G$ be a group.
An $G$-set $M$ is said to be \it connected \rm \cite{Bo97} if
there exists a connected $G$-graph $\Gamma$ with the vertex set $M$ such that
$\Gamma^1/G$ is finite. We call such $\Gamma$ a \it connecting structure \rm
for the $G$-set $M$.

Let $M$ be a connected $G$-set. A symmetric set $\b u{\subset}M^2$ containing $\b\Delta^2M$
is called a \it divider \rm(see \ref{expansiveMetric})
if there exists a finite set $F{\subset}G$ such that $(\cap\{f\b u:f{\in}F\})^2{\subset}\b u$.
The $G$-filter generated by a divider is a uniformity $\Cal U_{\b u}$ on $M$, see \ref{dividerEnt}.

A divider $\b u$ is said to be \it perspective \rm(see \ref{perspectivity})
if for every pair $\beta$ of points in $M$ the set $\{g{\in}G:g\beta{\notin}\b u\}$ is finite.

\deff{ $\s{RH_{pd}}$\rm. A \it relatively hyperbolic structure \rm
on a connected $G$-set $M$ (we also say
`a \it relatively hyperbolic uniformity\rm' or `a \it relative hyperbolicity\rm' on $M$)
is a $G$-uniformity $\Cal U$
generated by a perspective divider
such that the $\Cal U$-boundary has at least two points.}

There is no explicit expression of the ``parabolic'' subgroups for $\s{RH_{pd}}$.
We prove in \ref{pd2fh} that the completion $\overline{(M,\Cal U_{\b u})}$ with respect to a r.h. uniformity
is a compactum where $G$ acts relatively hyperbolically in the sense of $\s{RH_{32}}$.
So we obtain an interpretation $\s{RH_{pd}{\Rightarrow}RH_{32}}$.
% [map theorem] @tag_4DFCD6F4:
\subsection{Floyd map theorem for r.h. uniformities}
We use the exponential scaling function $f_n{=}\lambda^n$ where $0{<}\lambda{<}1$.

Let $\Gamma$ be a connected graph.
For a vertex $v{\in}\Gamma^0$ define the Floyd metric $\d_{v,\lambda}$
by postulating that the length of an edge of distance $n$ from $v$ is $\lambda^n$.
Change the base vertex $v$ gives
a bi-lipschitz-equivalent metric
and hence the same uniformity $\s U_{\Gamma,\lambda}$.
We call it the \it Floyd uniformity \rm on $\Gamma$.
\vskip4pt
\bf Theorem \rm(Map theorem \ref{mapThrm}) \sl
Let $G$ be a group, $M$ a connected $G$-set, $\Gamma$ a connecting graph structure for $M$,
$\Cal U$ a relatively hyperbolic uniformity on $M$.
Then there exists $\lambda\in(0,1)$ such that $\Cal U$
is contained in the Floyd uniformity $\s U_{\Gamma,\lambda}$.
The inclusion induces a uniformly continuous
$G$-equivariant surjective map
$\overline{(M,\s U_{\Gamma,\lambda})}\to\overline{(M,\Cal U)}$
between the completions\rm.
\vskip3pt
This theorem can be applied in particular to either a Cayley graph
with respect to a finite generating set or to a Farb's ``conned-off''
graph relative with respect to a finite collection of subgroups
without any restriction on the cardinality of the ``parabolic'' subgroups.
\vskip3pt
\bf Corollary\sl. Let $G$ be a group relatively hyperbolic with respect to
a collection $\Cal P$ of subgroups. Then, for some $\lambda{\in}(0,1)$, there exists
a continuous equivariant map from the Floyd boundary $\partial_\lambda G$ to the
Bowditch boundary of $G$ with respect to $\Cal P$\rm.
% [visibility] @tag_4E088B5C:
\subsection{Visibility}
For an edge $e$ of a graph $\Gamma$ let $\b u_e$ denote the set of pairs $(x,y)$
of vertices such that no geodesic segment joining $x$ and $y$ pass through $e$.
The filter $\s{Vis}\Gamma$ on the set of pairs of vertices
generated by the collection $\{\b u_e:e{\in}\Gamma^1\}$ is called the \it visibility filter\rm.
A graph $\Gamma$ is alternatively hyperbolic
if and only if $\s{Vis}\Gamma$ is a uniformity on $\Gamma^0$.

A uniformity $\Cal U$ on $\Gamma^0$ is called a \it visibility \rm on $\Gamma$
if it is contained in $\s{Vis}\Gamma$, see \ref{visibility}.
On every graph $\Gamma$ there exists the maximal visibility that contains all other visibilities.
It may be smaller than $\s{Vis}\Gamma$.

The main corollary of the map theorem is the following highly useful fact.
\vskip3pt
\bf Generalized Karlsson lemma \rm(\ref{KarlssonLemma})\sl.
Let $\Cal U$ be a relative hyperbolicity on a connected $G$-set $M$
and let $\Gamma$ be a connected graph with $\Gamma^0{=}M$ where $G$ acts on edges properly and cofinitely.
Then, for every entourage $\b u{\in}\Cal U$
there exists a finite set $E{\subset}\Gamma^1$ such that $\b u$
contains the boundary pair $\partial I$ of every geodesic segment $I$ that misses $E$\rm.
\vskip3pt
This implies that each relative hyperbolicity is a visibility.
\vskip3pt
The completion with respect to any visibility $\Cal U$ on a graph $\Gamma$
is compact (\ref{vizCompact}). If a group $G$ acts on $\Gamma$
properly on edges and keeps $\Cal U$ invariant
then the induces action on the completion space $T$ has the convergence property.
This gives a wide class of convergence group actions including the action on the space of ends,
Kleinian actions of finitely generated groups, the actions on Floyd completions and many other.
The problem is whether there exist convergence actions of other nature.

If $\Cal U$ is a relative hyperbolicity then $T$ coincides with the Bowditch completion,
see \ref{BoCmpl}.
In this case the induces action on the space of pairs is cocompact.
% [attractor sum] @tag_4E088B73:
\subsection{Attractor sum}
To prove `$\s{RH_{32}{\Rightarrow}RH_{pd}}$' we need to attach to a compactum $T$
where a group $G$ acts properly on triples at least one
orbit of isolated points.
We do so by a rather general construction that we call \it attractor sum\rm.

For the sake of future applications we construct this space in an excessive generality
of the actions of locally compact groups.
However the additional difficulties implied by possible non-discreteness of the acting group $G$
are not very essential. Moreover they clarify and motivate some aspects of the theory
of discrete group actions.

Our main result in this direction is the following.
\vskip3pt
\bf Attractor sum theorem \rm(\ref{attractorSum})\sl.
Let
a locally compact group $G$ act on a compactum $\Lambda$ properly on triples
and on a locally compact Hausdorff space $\Omega$ properly and cocompactly.
Then on the disjoint
union $\Lambda{\sqcup}\Omega$
there is a unique compact Hausdorff topology $\tau$
extending the original topologies of
$\Lambda$ and $\Omega$
such that the $G$-action on the space $X{\leftrightharpoons}(\Lambda{\sqcup}\Omega,\tau)$ is proper on triples\rm.
\vskip3pt
In particular every convergence group action of a discrete group $G$ on a compactum $T$
extends to a convergence group action on the attractor sum $\widetilde T$ of $G$ and $T$.
Thus the uniformity of $\widetilde T$ induces a uniformity on $G$.
One could ask whether this uniformity is a visibility or a relative hyperbolicity.
The closure $\overline G$ of $G$ in $\widetilde T$ can be thought of as an invariant compactification
of $G$ where $G$ acts with convergence property. The ``boundary'' $\overline G{\setminus}G$
is just the limit set of the original action $G{\on}T$.
% [structure] @tag_4E088B8C:
\subsection{The structure of the paper}\label{structPaper}
The following diagram presents the main results
illustrating the reasons and dependencies.

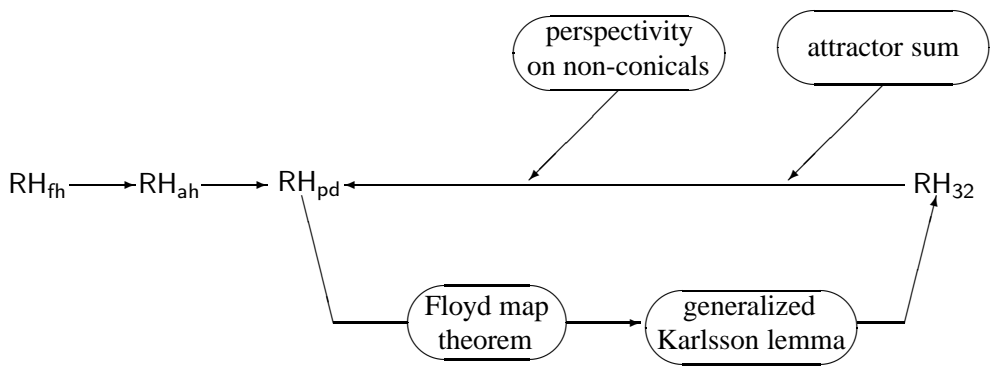
\begin{picture}(300,100)(-120,-40)
\put(120,50){\oval(80,28)}
\put(120,52){\makebox(0,0)[c]{\vbox{\hbox{\kern7pt perspectivity}\hbox{on non-conicals}}}}
\put(220,52){\oval(80,28)}
\put(220,52){\makebox(0,0)[c]{attractor sum}}
\put(120,36){\vector(-1,-1){34}}
\put(220,38){\vector(-1,-1){36}}
\put(-100,0){\makebox(0,0)[c]{$\s{RH_{fh}}$}}
\put(-88,0){\vector(1,0){25}}
\put(-50,0){\makebox(0,0)[c]{$\s{RH_{ah}}$}}
\put(-38,0){\vector(1,0){25}}
\put(3,0){\makebox(0,0)[c]{$\s{RH_{pd}}$}}
\put(228,0){\vector(-1,0){212}}
\put(243,0){\makebox(0,0)[c]{$\s{RH_{32}}$}}
\put(0,-4){\line(1,-4){12}}
\put(12,-52){\line(1,0){28}}
\put(70,-52){\oval(60,28)}
\put(70,-52){\makebox(0,0)[c]{\vbox{\hbox{Floyd map}\hbox{\kern5pt theorem}}}}
\put(100,-52){\vector(1,0){28}}
\put(170,-54){\oval(80,28)}
\put(170,-52){\makebox(0,0)[c]{\vbox{\hbox{\kern10pt generalized}\hbox{Karlsson lemma}}}}
\put(210,-52){\line(1,0){18}}
\put(228,-52){\vector(1,4){12}}
\end{picture}
\vskip33pt
The interpretation $\s{RH_{32}{\Rightarrow}RH_{fh}}$
for finitely generated groups follows from \cite{Ge09} and \cite{Ya04}.
The argument of \cite{Ya04} is rather complicated and it is not clear for us
whether or not the finite generability actually used.
Anyway it uses the metrisability of the compactum which is equivalent to the assumption
that the group is countable.

In \cite{GP10} we remove any restriction on the cardinality
and in \cite{GP11} we give a conceptually more simple proof of `$\s{RH_{32}{\Rightarrow}RH_{fh}}$'
using ``quasigeodesics'' with respect to a quadratic distortion functions.
In \cite{GP09} we study quasi-isometric maps to r.h. groups.

All the three papers require the map theorem \ref{mapThrm} and the attractor sum theorem \ref{attractorSum}
(for discrete groups).
% [introduction] @tag_4DEB8D8F) end of range
% [preliminaries] @tag_4A80F2BA(
\section{Preliminaries}\label{preliminaries}
In this section we fix the terminology and the notation and recall
widely known definitions and facts.
For the reader's convenience we repeat in this section the definitions given in the introduction.
The reader should search here for all general definition and notation used elsewhere in this paper.
We are trying to collect those and only those
definitions and statements that we use.
Some conventions are introduced implicitly.

We recommend to browse this section and to continue reading farther
returning to the preliminaries following the references.
% [notation] @tag_4A8237C3:
\subsection{General notation and conventions}
The symbol `$\square$' at the end of line means that the current proof
is either completed or left to the reader.
The reader is supposed to be capable to
complete the proof or to find it the common sources.

The single quotes `\dots' mean that the content is just mentioned, not used.
The double quotes ``\dots'' mean that the exact interpretation of the content is left to the reader.
Example: `$f_n$ tends to infinity' means ``$f_n$ gets arbitrarily big while $n$ grows''.

The symbol `$\leftrightharpoons$' means `is equal by definition'.
We use the \it italic font \rm for the notions being defined.
Example:
$\mathbb N{\leftrightharpoons}\mathbb Z_{\geqslant0}$ is the set of positive integers;
the elements of $\mathbb N$ are the \it natural numbers\rm.

For a set $M$ we denote by $M^n$ the product of $n$ copies of $M$.
The quotient of $M^n$ by the action of the symmetric group transposing the coordinates
is denoted by $\s S^nM$. The elements of $\s S^nM$ are the ``subsets of cardinality $n$ with multiplicity''.
The elements of $M^n$ and of $\s S^nM$ are the $n$\it-tuples\rm,
respectively ordered or non-ordered.
The 2-tuples are the \it pairs\rm, the 3-tuples are the \it triples \rm and so on.
An $n$-tuple is \it regular \rm if all its ``elements'' are distinct.
Denote $\b\Theta^nM\leftrightharpoons\{$regular $n$-tuples in $\s S^nM\}$,
$\b\Delta^nM\leftrightharpoons\s S^nM{\setminus}\b\Theta^nM=\{$\it singular \rm $n$-tuples\}.
We identify $\b\Theta^nM$ with $\{$subsets of $M$ of cardinality $n\}$.

By $|M|$ we denote the cardinality of a set $M$. However `$|M|{=}\infty$' means
'$M$ is infinite'.

For a set $M$ we denote by $\s{Sub}M$ the set of all subsets of $M$
and $\s{Sub}^nM{\leftrightharpoons}\{N{\subset}M:|N|{=}n\}$,
$\s{Sub}^{<n}M{\leftrightharpoons}\{N{\subset}M:|N|{<}n\}$ etc.
We identify $\s{Sub}^nM$ with $\b\Theta^nM$.

'$A{\setminus}B$' means set-theoretical difference.
For a subset $B$ of a set $A$ we write $B'$ instead of $A{\setminus}B$ when we hope that
the reader knows what is $A$.
Example: $B$ is open if and only if $B'$ is closed.

We sometimes identify a single-point set $\{\frak p\}$ with its unique point. For example,
for $\frak{p,q}{\in}T$ we write $\frak p{\times}T\cup T{\times}\frak q$ instead of
$\{\frak p\}{\times}T\cup T{\times}\{\frak q\}$. Speaking of the one-point compactification $\widehat L$
of a space $L$ we write $\widehat L{=}L{\cup}\infty$ instead of $\widehat L{=}L{\cup}\{\infty\}$ etc.

If $M$ is a Cartesian product of sets $M_\xi$ then a \it subproduct \rm of $M$ is a subset
which is a product of a family $N_\xi$ of subsets. A subproduct of the form $\pi_\xi^{-1}N$
where $\pi_\xi:M\to M_\xi$ is the projection map, is the \it cylinder over \rm$N$. If $N{=}\{\frak p\}$
then it is the \it fiber over \rm the point $\frak p$.

We use $\mathbb{R,Z}$ etc in the common way.
For $a,b{\in}\mathbb R$ by $[a,b]$, $[a,b)$, $(a,b)$ we denote
the intervals of different types, as the reader expects.
For $m,n{\in}\mathbb Z$ we put $\overline{m,n}\leftrightharpoons\{k{\in}\mathbb Z:m{\leqslant}k{\leqslant}n\}=[m,n]{\cap}\mathbb Z$.

By $f|_M$ we denote the restriction of a function $f$ onto a set $M$ not necessarily contained in
the domain $\s{Dom}(f)$ of $f$.
So, $\s{Dom}(f|_M){=}M{\cap}\s{Dom}(f)$.

We consider an equivalence relation on a set $M$ as a subset of either $M^2$ or $\s S^2M$.
The same convention is adopted for other ``symmetric'' relations and ``symmetric'' functions.
The reader should not be confused.

The \it kernel \rm of a map $f$ defined on a set $M$ is the equivalence relation `$f\frak p{=}f\frak q$'.

For a function $f$ on a subset of $\mathbb Z$ we sometimes write $f_n$ instead of $f(n)$.
% [topology] @tag_4A9D2B8A:
\subsection{General topology}\label{gTopology}
The information of this subsection is used mainly in the attractor sum theory of section \ref{attSumSect}.

For a topological space $T$ we denote:
$\s{Open}(T){\leftrightharpoons}\{$open subsets of $T\}$;\hfil\penalty-10000
$\s{Closed}(T){\leftrightharpoons}\{$closed subsets of $T\}$;
$\s{Loc}_TS{\leftrightharpoons}\{$the neighborhoods of a subset $S$ of $T\}$.
A \it filter \rm on a set $M$ is a \bf proper \rm subset $\Cal F$ of $\s{Sub}M$
such that $A,B{\in}\Cal F{\Leftrightarrow}A{\cap}B{\in}\Cal F$.
An \it ultrafilter \rm on $M$ is a maximal element in the set \{filters on $M\}$.
Every filter is the intersection of a family of ultrafilters.
A collection of sets is \it consistent \rm if any finite subcollection
has nonempty intersection.
A collection is consistent if and only if it is contained in some filter.

Two collections are \it inconsistent \rm if their union is not consistent.

A filter on a topological space $T$ \it converges \rm to a point $\goth p{\in}T$ if it contains
$\s{Loc}_T\goth p$. A space $T$ is \it compact \rm if every
ultrafilter on $T$ converges (to a point). A Hausdorff compact space is a \it compactum\rm.

A subset of a topological space is (topologically) \it bounded \rm if its closure is compact.

\begin{prop}\label{closedGenFil}
If a filter $\mathcal F$ on a compact space $T$ is generated by $\mathcal F{\cap}\s{Closed}T$
then $\mathcal F{\supset}\s{Loc}_T({\cap}\mathcal F)$.\end{prop}\sl Proof\rm.
Every ultrafilter containing $\Cal F$ converges to a point that must belong to ${\cap}\Cal F$.\eod
\begin{prop}\label{generatingLoc}Let a space $T$ be compact, $F{\in}\s{Closed}T$
and let $\Cal F$ be a subfilter of $\s{Loc}_TF$ inconsistent with each $\s{Loc}_T\goth p$, $\goth p{\notin}F$.
Then $\Cal F{=}\s{Loc}_TF$.\end{prop}\sl Proof\rm.
Every ultrafilter that contains $\Cal F$ must converge to a point in $F$.\eod
\begin{prop}\label{ultra}
If an ultrafilter $\Cal F$ on a space $T$ contains $\s{Loc}_T\Lambda$ for a compact
subspace $\Lambda$ then
$\Cal F$ converges to a point of $\Lambda$.\end{prop}\sl Proof\rm.
The set $\Cal G{\leftrightharpoons}\Cal F{\cap}\s{Closed}T$ is consistent with $\{\Lambda\}$
thus there exists $\goth p{\in}\Lambda{\cap}({\cap}\Cal G)$.\hfil\penalty-10000
Every $o{\in}\s{Open}T{\cap}\s{Loc}_T\goth p$ should belong to $\Cal F$ since its complement
does not belong to $\Cal G$.\eod
\vskip3pt
The \it product topology \rm on a product $X$ of topological spaces is generated
by the set \{cylinders over open sets\}.
\begin{prop}[``Walles Theorem'']\label{Walles}
If $T$ is a compact subproduct of a product $X$ of topological spaces
then the filter $\s{Loc}_XT$
is generated by $\mathcal F{\cap}\{$cylinders$\}$.\eod\end{prop}
\begin{prop}[``Aleksandrov Theorem'']\label{Aleksandrov} The quotient $X/\theta$
of a compact space $X$ by an equivalence $\theta$ is Hausdorff
if and only if $\theta$ is closed in $X^2$.\eod\end{prop}
\begin{prop}[``Kuratowski Theorem'']\label{kuratowski} A topological space $K$ is compact if and only if
for every space $T$ the projection map $T{\times}K\to T$ is \it closed \sl i.e, maps closed sets to closed sets.\eod\end{prop}
A \it closed correspondence \rm from a space $X$ to a space $Y$ is any closed subset $S$ of $X{\times}Y$.
It is \it surjective \rm if the restrictions over $S$ of the projections onto $X$ and $Y$ are surjectve.

Let $A,B$ be topological spaces and let $K$ be a compactum.
For $a{\in}\s{Closed}(A{\times}K)$, $b{\in}\s{Closed}(K{\times}B)$ the set
$a{*}b\leftrightharpoons a{\times}b\cap A{\times}\b\Delta^2K{\times}B$
is closed in $A{\times}K^2{\times}B$.
The \it composition \rm$a{\circ}b{\leftrightharpoons}\s{pr}_{A\times B}a{*}b$ is closed by \ref{kuratowski}.

Let $\s{Surj}(X{\times}Y){\leftrightharpoons}\{$surjective closed correspondences from $X$ to $Y\}$.
\begin{prop}\label{composition}The composition of surjective correspondences is surjective.
The operation
$\s{Surj}(A{\times}Y)\times\s{Surj}(Y{\times}B)\to\s{Surj}(A{\times}B)$
is associative.\eod\end{prop}
% [metrics] @tag_4A835192:
\subsection{Metrics and uniformities}\label{uniform}
\setcounter{equation}0
We extend the common notion of a metric by allowing infinite distance and zero distance between
different points.
So, a \it metric \rm on a set $M$ is a function
$\varrho:\s S^2M\to\mathbb R_{\geqslant0}{\cup}\infty$ with $\varrho|_{\b\Delta^2M}{=}0$ satisfying
the $\triangle$\it-inequality \rm$\varrho(a,b){+}\varrho(b,c)\geqslant\varrho(a,c)$.
A metric is \it finite \rm if it does not take the value $\infty$. A metric $\varrho$ on $M$
is \it exact \rm
if $\varrho^{-1}0{=}\b\Delta^2M$.

Let $\rho$ be a metric on $M$. We extend the function $\rho$ to the pairs of subsets of $M$:
$\rho(A,B){\leftrightharpoons}\s{inf}\rho|_{A{\times}B}$.
The \it open \rm and \it closed $r$-neighborhoods \rm of a set $A{\subset}M$
are the sets $\s N_\rho(A,r){\leftrightharpoons}\{\frak p{\in}M:\rho(A,\frak p){<}r\}$
and $\overline{\s N}_\rho(A,r){\leftrightharpoons}\{\frak p{\in}M:\rho(A,\frak p){\leqslant}r\}$.
Sometimes we omit the index `$\rho$'.

Any set $\b u{\subset}\s S^2M$ can be thought as a symmetric binary relation on $M$
and as a set of the edges of a graph whose vertex set is $M$.
We call $\b u$ \it reflexive \rm if $\b u{\supset}\b\Delta^2M$.
Every $\b u{\subset}\s S^2M$ determines a metric $\d_{\b u}$ on $M$ as the maximal among the metrics
$\varrho$ on $M$ such that $\varrho|_{\b u}\leqslant1$. The canonical graph metric $\s d$
discussed in \ref{grphs}
is a particular case of this construction.
Denote $\b u^n{\leftrightharpoons}\d_{\b u}^{-1}(\overline{0,n})$.
Clearly\begin{equation}\label{capExp}(\b{u{\cap}v})^n\subset\b u^n{\cap}\b v^n\text{ for every }\b{u,v}{\subset}\s S^2M.\end{equation}

A set $m{\subset}M$ is $\b u$\it-small \rm if $\s S^2m\subset\b u{\cup}\b\Delta^2M$ (equivalently: if its $\d_{\b u}$-diameter
is ${\leqslant}1$). We denote by $\s{Small}(\b u)$ the set of $\b u$-small subsets of $M$.
We try to use the convention ``small letters denote small sets''.
 The $\b u$\it-neighborhood \rm of a set $m{\subset}M$ is the set
$m\b u{\leftrightharpoons}\overline{\s N}_{\d_{\b u}}(m,1)=m{\cup}\{\frak p{\in}M:\exists\frak q{\in}m\ \{\frak{p,q}\}{\in}\b u\}$.

Subsets $\b{u,v}{\subset}\s S^2M$ are said to be \it unlinked \rm if $M$ is a union of a $\b u$-small set and a $\b v$-small set.
We denote this relation by $\b{u{\bowtie}v}$. If $\b u$ and $\b v$ are not unlinked we say that they are \it linked \rm
and denote this relation by $\b{u{\#}v}$. So, $\b u$ is \it self-linked \rm($\b{u{\#}u}$) if $M$ is not a union of two
$\b u$-small sets.
% [uniformity] @tag_4A8351AA:
A filter $\mathcal U$ consisting of reflexive subsets of $\s S^2M$ is
a \it uniformity \rm or a \it uniform structure \rm on $M$
if $\forall\b u{\in}\mathcal U\exists\b v{\in}\mathcal U\ \b v^2{\subset}\b u$.

The elements of a uniformity are called \it entourages\rm.
This notion plays a significant role in our theory.
We use the \bf bold \rm font for the entourages and for some sets of pairs that ``should be'' entourages of some uniformities
(see for example, subsections \ref{Frink}, \ref{expansiveMetric}).

If $\b u$ is an entourage of a uniformity $\Cal U$ we write $\b v{=}\root n\of{\b u}$
if $\b v{\in}\Cal U$ and $\b v^n{\subset}\b u$. So $\root n\of{\b u}$ exists but it is not unique.

An entourage $\b u$ \it separates points $x$ and $y$ if $\{x,y\}{\notin}\b u$.
A uniformity $\mathcal U$ on $M$ is \it exact \rm if
every two distinct points can be separated by an entourage, i.e, ${\cap}\mathcal U{=}\b\Delta^2M$.

Given a uniformity $\mathcal U$
a set $m{\subset}M$ is called a $\mathcal U$\it-neighborhood \rm of a point $\frak p{\in}M$ if it contains
a $\b u$-neighborhood $\frak p\b u$ for some entourage $\b u{\in}\mathcal U$.
So $\mathcal U$ yields the $\mathcal U$\it-topology \rm on $M$ in which the neighborhoods of points are
the $\mathcal U$-neighborhoods. We speak of $\mathcal U$\it-open \rm sets $\mathcal U$\it-closure \rm etc. meaning the
$\mathcal U$-topology.

A topological space whose topology is determined by a uniformity $\mathcal U$ is \it uniformisable\rm.
Every such $\mathcal U$ is a uniformity \it consistent \rm with the topology.

For every compactum (moreover, even for every paracompact Hausdorff space) $T$ the filter\hfil\penalty-10000
$\s{Ent}(T){\leftrightharpoons}\s{Loc}_{\s S^2T}\b\Delta^2T$
of the neighborhoods of the diagonal is an exact uniformity on $T$ consistent with the topology.
If $T$ is compactum then the uniformity consistent with the topology is unique.
Therefore we have a correct notion of an entourage of a compactum $T$. It is just
a neighborhood of the diagonal $\bold\Delta^2T$ in the space $\mathsf S^2T$.

A topological space is uniformisable by an exact uniformity if and only if it is embeddable in a compactum
\cite[$\mathchar"278$1 Prop. 3]{Bou58}.

A set endowed with a uniformity is a \it uniform space \rm\cite{Bou71}, \cite{We38}.

Every metric $\varrho$ on $M$ determines a uniformity $\mathcal U_\varrho$ generated by the collection
$\{\varrho^{-1}[0,\varepsilon]:\varepsilon{>}0\}$. A uniformity is determined by a metric (=\it metrisable\rm)
if and only if it is countably generated as a filter (\cite{We38}, see also \ref{Frink})

A metric is exact if and only if the corresponding uniformity is exact.

The morphisms of uniformities are the \it uniformely continuous maps \rm i.e, the maps
such that the preimage of an entourage is an entourage. Every subset $N$ of a uniform space $(M,\Cal U)$
has the induces uniformity $\Cal U|_N$. It is the minimal among the uniformities on $N$ for which
the inclusion is uniformly continuous. The space $(N,\Cal U|_N)$ is a \it subspace \rm of $(M,\Cal U)$.
% [completion] @tag_4DECB2F1:
\subsection{Cauchy-Samuel completion}\label{Samuel}\setcounter{equation}0
Let $(M,\Cal U)$ be a uniform space. A \it Cauchy filter \rm$\Cal F$ on $M$ is a filter
with arbitrarily small elements:
 $\forall\b u{\in}\Cal U\,\Cal F{\cap}\s{Small}(\b u){\ne}\varnothing$.
For $x{\in}M$ the filter $\s{Loc}_{\Cal U}x$ is a Cauchy filter.
Moreover, it is minimal element in the set of Cauchy filters ordered by inclusion.

The space is \it complete \rm if every Cauchy filter \it converges \rm i.e, contains a filter
of the form $\s{Loc}_{\Cal U}x$ for $x{\in}M$.
Every closed subset of a complete space is a complete subspace.
Every compactum is complete.

Every uniform space $(M,\Cal U)$ possesses an initial morphism
$\iota_{\Cal U}:(M,\Cal U)\to(\overline M,\overline{\Cal U})$
to a complete space. The points of $\overline M$ are the minimal Cauchy filters.
The completion map $\iota_{\Cal U}$ takes $x$ to $\s{Loc}_{\Cal U}x$.
For an entourage $\b u{\in}\Cal U$ the set
\begin{equation}\label{Samuel1}
\overline{\b u}{\leftrightharpoons}\{\g{p,q}\}{\in}\s S^2\overline M:
\g{p{\cap}q}{\cap}\s{Small}(\b u){\ne}\varnothing\}\end{equation}
is, by definition, an entourage of $\overline M$.
The uniformity $\overline{\Cal U}$ is the filter generated by
$\{\overline{\b u}:\b u{\in}\Cal U\}$. It is exact:
if $\g{p{\ne}q}$ then the filter
$\g{p{\cap}q}$ is smaller than some of $\g{p,q}$
since it is not a Cauchy filter. So $\{\g{p,q}\}$ does not belong to some $\overline{\b u}$.

If $\Cal U$ is exact then $\iota_{\Cal U}$ is injective
and we can identify $M$ with a subspace of $\overline M$.
In this case the \it remainder \rm$\partial_{\Cal U}M{\leftrightharpoons}\overline M{\setminus}M$
is sometimes called a $\Cal U$\it-boundary \rm of $M$.

For every subset of a complete exact uniform space the canonical map from the completion
to the closure is an isomorphism of uniform spaces.

An entourage $\b u$ is \it precompact \rm if $M$ is a union of finitely many $\b u$-small sets.
A uniformity is \it precompact \rm if every its entourage is.
Precompactness of a uniformity is equivalent to the compactness of the completion space
\cite[Thm 32, p. 198]{Ke75}.
\subsection{Dynkin property}
Two entourages $\b{u,v}$ of a uniform space $(M,\Cal U)$ are said to be \it unlinked \rm
(notation: $\b{u{\bowtie}v}$)
if $M$ is a union of an $\b u$-small set and a $\b v$-small set. Otherwise
the entourages are \it linked \rm(notation: $\b{u{\#}v}$).

Let a locally compact group $G$ act on $M$ keeping $\Cal U$ invariant.
We say that the action $G{\on}M$ \it has Dynkin property \rm if
for every $\b{u,v}{\in}\Cal U$ the set $\{g{\in}G:\b u{\#}g\b v\}$ is bounded in $G$.
For compact spaces the Dynkin property is equivalent to the ``convergence property'',
see \ref{charProp}.
% [complDynkin] @tag_4DF3AA00
\begin{prop}\label{complDynkin}Completion keeps Dynkin property.\end{prop}\sl Proof\rm.
It suffices to check that,
for unlinked entourages $\b{u,v}$ of a uniform space $(M,\Cal U)$,
the entourages $\overline{\b u^3}$ and $\overline{\b v^3}$ (see \ref{Samuel1}),
are unlinked.

For a set $a{\subset}M$ let $\widetilde a{\leftrightharpoons}\{\g p{\in}\overline M:a$ is consistent with $\g p\}$
If $a$ is $\b u$-small then $\widetilde a$ is $\b u^3$-small.
Indeed if $\g{p,q}\in\widetilde a$ and $p,q$ are $\b u$-small sets in $\g{p,q}$ respectively
then the set $p{\cup}a{\cup}q$ is an $\b u^3$-small set in $\g{p{\cap}q}$.

If $\b{u{\bowtie}v}$ and $M{=}a{\cup}b$ where $a{\in}\s{Small}(\b u)$, $b{\in}\s{Small}(\b v)$ then
every filter on $M$ is consistent with either $a$ or $b$.
So $\overline M{=}\widetilde a{\cup}\widetilde b$.\eod

% [graphs] @tag_4DFB08D5:
\subsection{Graphs}\label{grphs} For a graph $\Gamma$ we denote by $\Gamma^0$ and $\Gamma^1$ the sets of vertices
and edges respectively. We do not interesting in graphs with loops and multiple edges.
For our purpose a graph is something that is either connected or disconnected.
We often identify an edge $e{\in}\Gamma^1$ with its boundary pair $\partial e{\subset}\Gamma^0$
and write ${\cup}E$ for the set ${\cup}\{\partial e:e{\in}E\}$ ($E{\subset}\Gamma^1$).
We consider any set of pairs as a graph.

By $\s d$ or by $\s d_\Gamma$ we denote the natural metric on $\Gamma^0$ it is the maximal
among the metrics for which the distance between joined vertices is one.

A \it circuit \rm is a connected graph with exactly two edges at each vertex.
An \it arc \rm is a graph obtained from a circuit by removing one edge.
By $\s{Arc}(\Gamma,e)$ we denote the set of all arcs in $\Gamma$ that contain the edge $e{\in}\Gamma^1$.

A graph $\Gamma$ is \it fine \rm\cite{Bo97} if each set of arcs of bounded length with fixed endpoints is finite.
% [perspectivity] @tag_4DFA742C:
\subsection{Perspectivity}\label{perspectivity}
When something goes farther it looks smaller. This phenomenon is the \it perspectivity\rm.
\deff{. A uniformity $\Cal U$ on the set of vertices of a connected graph $\Gamma$ is
said to be \it perspective \rm
if every entourage contains all but finitely many edges. We also say that $\Cal U$ is a \it perspectivity \rm on $\Gamma$.}
\begin{prop}\label{exactIsFine}
If a graph $\Gamma$ possesses an exact perspectivity $\Cal U$ then it is fine.\end{prop}\sl Proof\rm.
Suppose that $\Gamma$ is not fine.
Let $n$ be the smallest positive integer for which
there exists an infinite set $P$ of arcs of length at most $n$ joining
some fixed vertices $x,y{\in}\Gamma^0$.
Let an entourage $\b u$ separate $\{x,y\}$ and let $\b v{=}\root n\of{\b u}$.
We have $P\subset{\cup}\{\s{Arc}(\Gamma,e):e{\in}\Gamma^1{\setminus}\b v\}$.
By the perspectivity property the set $\Gamma^1{\setminus}\b v$ is finite.
So, for some fixed $e{\in}\Gamma^1{\setminus}\b v$, the set $P{\cap}\s{Arc}(\Gamma,e)$
is infinite that gives us a counter-example
with a smaller $n$.\eod
\vskip3pt
Let a locally compact group $G$ act on a uniform space $(M,\Cal U)$.
An (ordered or non-ordered) pair $\beta$ of points of $M$ is
said to be \it perspective \rm if the orbit map $G{\ni}g\mapsto g\beta{\in}\b\Theta^2M$
is ``proper'', i.e, for every entourage $\b u{\in}\Cal U$ the set $\{g{\in}G:g\beta{\notin}\b u\}$
is bounded in $G$.

Perspectivity is a property of an \bf orbit \rm in $\b\Theta^2M$.
On the other hand it is an equivalence relation.
Hence if $M{=}\Gamma^0$ is the vertex set of a \bf connected \rm graph $\Gamma$
and $G$ acts by graph automorphisms then the perspectivity on edges implies the perspectivity of the action.

On the other hand perspectivity is a relation between a pair and an orbit of entourages.
If a $G$-uniformity $\Cal U$ is generated by a set $\Cal S$ of entourages
and $\beta$ is perspective with respect to each $\b u{\in}\Cal S$ then $\beta$ is perspective with respect to $\Cal U$.
\begin{prop}\label{properOnPairs}
Let $\Cal U$ be an invariant exact perspectivity on a connected $G$-graph $\Gamma$ where $G$ acts properly on edges.
Then $G$ acts properly on pairs.\end{prop}\sl Proof\rm.
By \ref{exactIsFine} $\Gamma$ is fine so if $\{x,y\}{\in}\b\Theta^2\Gamma^0$
then the stabilizer $\s{St}_G\{x,y\}$ acts properly on the finite set
of the geodesic arcs between $x$ and $y$.\eod
\begin{prop}\label{perspectivityCompl}
Let $\Cal U$ be an invariant perspectivity on a connected $G$-graph $\Gamma$ where $G$ acts properly on edges.
Let $\iota_{\Cal U}:\Gamma^0\to\overline{\Gamma^0}$ be the completion map and
let $\Delta^0{\leftrightharpoons}\iota_{\Cal U}\Gamma^0$ and
$\Delta^1{\leftrightharpoons}(\iota_{\Cal U}\Gamma^1){\setminus}\{$loops\}.
The action $G{\on}\Delta$ is proper on pairs.\end{prop}\sl Proof\rm.
\note{I am grateful to the referee suggested the following elegant proof}
By \ref{properOnPairs} it suffices to prove that the action is proper on edges.
Let $e{=}\{p,q\}{\in}\Gamma^1$. Suppose that the filters $\iota_{\Cal U}p{\leftrightharpoons}\g p$,
$\iota_{\Cal U}q{\leftrightharpoons}\g q$
are distinct.
If $g\{\g{p,q}\}\notin\overline{\b u}$
then, by \ref{Samuel1},
$g\g p{\cap}g\g q{\cap}\s{Small}(\b u){=}\varnothing$.
Let $\b v{=}\root3\of{\b u}$.
The pair $\{gp,gq\}$ is not $\b v$-small
and $g$ belongs to a bounded subset of $G$.\qed
% [changing generating set] @tag_4A80F27F:
% [preliminaries] @tag_4A80F2BA) end of range
% [Floyd map] @tag_4A80F403(
\section{Floyd map}
% [Frink] @tag_4A813502:
\subsection{Frink Lemma}\label{Frink}
\setcounter{equation}0
The following lemma is a well-known metrisation tool of the general topology
(\cite[$\mathchar"278$1 Prop. 2]{Bou58}, \cite[Lemma 6.2]{Ke75}, \cite[Theorem 8.1.10]{En89}).
J. L. Kelly attributes it to A.~H.~Frink \cite{Fr37} noting certain contribution of
other authors. R. Engelking cites a paper \cite{Tu40} of another author.

A \it Frink sequence \rm on a set $M$ is a sequence $\b v_n$ such that $\b v_0{=}\s S^2M$
and $\b v_n{\supset}\b v_{n+1}^3{\supset}\b\Delta^2M$ for all $n{\in}\mathbb N$.
Any Frink sequence $\b v_*$ determines on $M$ the \it Frink metric \rm
as the maximal among the metrics $\varrho$ on $M$ such that
$\forall n\ \varrho|_{{\fam6v}_n}\leqslant2^{-n}$ (recall that we do not require for a metric to be exact).
\begin{prop}[Frink Lemma]\label{frinkLemma} Let $\d$ denote the Frink metric on a set $M$
determined by a Frink sequence $\b v_*$.
Then $\forall n{>}0\ \d^{-1}[0,2^{-n})\subset\b v_{n-1}$.
In particular the filter generated by $\b v_*$ is the uniformity determined by $\d$.
\end{prop}
For the reader's convenience we adopt the common proof to our notation.
The following definitions are valid only within the proof.

An \it edge \rm of a set $F{\subset}\mathbb R$
is a pair $(x,y){\in}F^2$ such that $x{<}y$ and $F{\cap}(x,y){=}\varnothing$.
A \it path \rm is a map
\setcounter{equation}1%
\begin{equation}\label{path}\gamma:F\to M\end{equation} from a \bf finite \rm$F{\subset}\mathbb R$
such that
$y{-}x=2^{-{\s{max}}\{n:\gamma\{x,y\}{\in}\b v_n\}}$ for every edge $(x,y)$ of $F$.
Write\hfil\penalty-10000
$\partial\gamma{\leftrightharpoons}\gamma\{\s{min}F,\s{max}F\}$, $\s{length}(\gamma){\leftrightharpoons}\s{max}F{-}\s{min}F$.
\vskip5pt
\sc Lemma\sl. If $\s{length}(\gamma)<2^{-n}$ for a path \ref{path} then
$\partial\gamma\in\b v_{n-1}$.
\par\vskip5pt\noindent\sl Proof\rm.
Induction by $|F|$. If $|F|{=}1$ then trivially $\partial\gamma{\in}\b\Delta^2M{\subset}\b v_{n-1}$.
Consider a path \ref{path} with\hfil\penalty-10000
$0{<}\s{length}(\gamma){\leftrightharpoons}l{<}2^{-n}$.
Suppose that the assertion is true for the proper subpaths of $\gamma$.

Let $(f_-,f_+)$ be the edge of $F$, for which ${1\over2}(\s{max}F{+}\s{min}F)\in[f_-,f_+]$.
The length of the restrictions of $\gamma$ over the sets $F{\cap}\mathbb R_{\leqslant f_-}$ and
$F{\cap}\mathbb R_{\geqslant f_+}$ is at most ${l\over2}{<}2^{-n-1}$.
By the inductive hypotheses
$\gamma\{\s{min}F,f_-\},\gamma\{f_+,\s{max}F\}\in\b v_n$.
Furthermore, $f_+{-}f_-\leqslant l{<}2^{-n}$, hence, by definition of path, $\gamma\{f_-,f_+\}{\in}\b v_n$.
So, by the definition of Frink sequence,
$\partial\gamma{\in}\b v_n^3{\subset}\b v_{n-1}$.\eod
\vskip5pt
The Frink distance between $\frak{p,q}{\in}M$ is, clearly,
the infimum of $\s{length}(\gamma)$ over all paths \ref{path} such that $\partial\gamma{=}\{\frak{p,q}\}$.
The assertion is now follows from the Lemma.\eod
% [divider] @tag_4DFA31D1:
\subsection{Frink sequence determined by a divider}\label{expansiveMetric}
\setcounter{equation}0
For a $G$-set $M$ a set $\b u{\subset}\s S^2M$ and a finite set $F{\subset}G$
denote $F\{\b u\}{\leftrightharpoons}\{f\b u:f{\in}F\}$. So ${\cap}(F\{\b u\})$ is a subset of $\s S^2M$.
It follows from \ref{capExp} that\begin{equation}\label{capExp1}({\cap}(F\{\b u\}))^n\subset{\cap}(F\{\b u^n\})\text{
for }n{\geqslant}1.\end{equation}
\deff{ \ref{expansiveMetric}.
We call a set $\b u{\subset}\s S^2M$ \it divider\rm, if it contains the diagonal $\b\Delta^2M$
and there exists
a finite $F{\subset}G$ such that $({\cap}(F\{\b u\}))^2{\subset}\b u$.}

Example: every equivalence relation is a divider for $F{=}\{1_G\}$.
We will see that every r.h. group contains a divider that determines
the relatively hyperbolic structure.
\setcounter{prop}1%
\begin{prop}\label{dividerEnt} Let $\b u$ be a divider for a connected $G$-set $M$.
The $G$-filter $\Cal U$ on $\s S^2M$ generated
as a $G$-filter by a set $\{\b u\}$ is a uniformity.\end{prop}\sl Proof\rm.
The sets of the form ${\cap}(S\{\b u\})$, $S{\in}\s{Sub}^{<\infty}G$ is a base for the filter
$\Cal U$. It follows from \ref{capExp1} that
$({\cap}(S{\cdot}F\{\b u\}))^2{\subset}\cap(S\{\b u\})$.\eod
\deff{ 3.2.3.
We say that $\Cal U$ is the uniformity \it generated \rm by $\b u$. We denote it by $\Cal U_{\b u}$.}

Every divider $\b u$ satisfies a stronger condition
\setcounter{equation}3%
\begin{equation}\label{expansive}
\forall m{>}0\exists F{\subset}G\ |F|{<}\infty\ \land\ ({\cap}(F\{\b u\}))^m{\subset}\b u.
\end{equation}
Indeed, by iterating the inclusion $({\cap}(F\{\b u\}))^2{\subset}\b u$ we obtain \ref{expansive} for $m{=}2^k$.
On the other hand $({\cap}(F\{\b u\}))^m\subset({\cap}(F\{\b u\}))^n$ for $m{<}n$ since $\b u$ contains
the diagonal.\eod

Note that if \ref{expansive} holds for a fixed $m$ and $F$ then
it holds for the same $m$ and every finite $F_1{\supset}F$.

Let $\b u$ be a divider
and let $F$ be a finite symmetric (i.e. closed under $g\mapsto g^{-1}$)
subset of $G$ containing the neutral element $1_G$ and such that
\begin{equation}\label{idemp3}
({\cap}(F\{\b u\}))^3{\subset}\b u.
\end{equation}
For $n{\in}\mathbb N$ let $F^n{\leftrightharpoons}\{$the elements of the group generated by $F$ of $F$-length ${\leqslant}n\}$.
The sets $F^n$ are symmetric.
It follows from the associativity that
\begin{equation}\label{expF} F^n{\cdot}F^m{=}F^{n+m}\text{ for all }m,n{\geqslant}0.
\end{equation}
The sequence\begin{equation}\label{divSeq}
\b u_0=\s S^2M,\ \b u_n\leftrightharpoons{\cap}(F^{n-1}\{\b u\})\text{ for }n{>}0\end{equation}
is a Frink sequence (see \ref{Frink}) containing the divider $\b{u{=}u}_1$.
Indeed, `$\b u_1^3{\subset}\b u_0$' is trivial and `$\b u_2^3{\subset}\b u_1$' is \ref{idemp3}.
Further, \ref{expF} implies
\begin{equation}\label{expU}
\b u_{n+m}={\cap}(F^n\{\b u_m\})\text{ for all }n{\geqslant}0, m{\geqslant}1.
\end{equation}
For $n{\geqslant}3$, using \ref{expU} and \ref{capExp1} we obtain
 $$\b u_n^3=({\cap}((F^{n-2}\{\b u_2\}))^3\subset{\cap}(F^{n-2}\{\b u_2^3\})\subset{\cap}(F^{n-2}\{\b u\})=\b u_{n-1}.$$
The Frink metric determined by this sequence is called the \it dividing metric\rm. It depends on $\b u$ and $F$.
The corresponding uniformity is called the \it dividing uniformity\rm.
% [comparison fl/fr] @tag_4A811611:
\subsection{Comparing the Floyd metric with the dividing metric}\label{cmpMetrics}\hfil\penalty-10000
\vskip-5pt
\deff{ \ref{cmpMetrics}.
A divider $\b u$ on a connected $G$-set $M$ is said to be \it perspective \rm
(compare with \ref{properOnPairs}) if
for each pair $\beta{\subset}M$ it contains all but finitely many elements of the orbit $G\{\beta\}$.}
Suppose that $M$ is the vertex set $\Gamma^0$ of a connected graph $\Gamma$
where $G$ acts properly on edges and cofinitely.
A divider $\b u$ is perspective if and only if it contains all but finitely many $\Gamma$-edges.

We fix $\Gamma$ and a perspective divider $\b u$.
Let a finite set $F$ satisfy \ref{idemp3} and denote by $\d_{\b u,F}$ the dividing metric.

 Denote by $\s d$ the canonical graph metric on $M$
(the ``$\s d$-length'' of each edge is 1).
For a vertex $v{\in}M$ and $\lambda{\in}(0,1)$ denote by $\d_{v,\lambda}$ the path metric on $M$ for which
the ``length'' of an edge $e$ is $\lambda^{\s d(v,e)}$.

\begin{prop}\label{cmpMetr} There exist $\lambda{\in}(0,1)$ and $C{>}0$ such that
$\d_{\b u,F}{\leqslant}C\d_{v,\lambda}$ on $M$.\end{prop}\sl Proof\rm.
Let $B_n{\leftrightharpoons}\{w{\in}M:\s d(v,w){\leqslant}n\}$.
The finite set $F\{v\}$ is contained in some $B_\rho$.
This implies $FB_n{\subset}B_{n+\rho}$ for all $n{\geqslant}0$.
Hence\setcounter{equation}1%
\begin{equation}\label{balls}F^kB_n{\subset}B_{n+k\rho}\text{ for }k,n{\geqslant}0.\end{equation}

Let $\sigma$ be such that $B_\sigma$ contains all the edges ${\notin}\b u$.
We put\begin{equation}\label{floydConstants}\lambda\leftrightharpoons2^{-1/\rho},\
C\leftrightharpoons2^{\sigma/\rho}.\end{equation}
Let $g{\in}F^k$. Since $g^{-1}{\in}F^k$, \ref{balls} implies
$\varnothing{=}B_\sigma{\cap}g(\Gamma^0{\setminus}B_{\sigma+k\rho})$ and hence
$g(\Gamma^1{\setminus}\s S^2B_{\sigma+k\rho}){\subset}\b u$.
We thus have\begin{equation}\label{edges}\Gamma^1{\setminus}\s S^2B_{\sigma+k\rho}\subset
{\cap}(F^k\{\b u\}){=}\b u_{k+1}.\end{equation}

We now compare the $\d_{\b u,F}$-length of an edge $e{\in}\Gamma^1$ with its $\d_{v,\lambda}$-length
for $\lambda{\leftrightharpoons}2^{-1/\rho}$.
Let $e$ be an edge of $\Gamma$ not contained in $B_\sigma$.
Then there exists a unique number $k{\geqslant}0$ such that $e\subset B_{\sigma+(k+1)\rho}{\setminus}B_{\sigma+k\rho}$.
By \ref{edges}  and the definition \ref{Frink} of the Frink length,
we have $\d_{\b u,F}(e){\leqslant}2^{-k-1}$.
The Floyd length of $e$ is $\d_{v,\lambda}(e)\geqslant\lambda^{\sigma+(k+1)\rho}=2^{-\sigma/\rho}{\cdot}2^{-k-1}$.
So we have\begin{equation}\label{FrFl}\d_{\b u,F}(e)\leqslant2^{\sigma/\rho}{\cdot}\d_{v,\lambda}(e)\end{equation}
for every edge outside $B_\sigma$.
The edges inside $B_\sigma$ also satisfy \ref{FrFl} since $\d_{\b u,F}{\leqslant}1$
everywhere.

Since the Floyd metric is maximal among those having the given length of edges,
and $\Gamma$ is connected,
\ref{FrFl} holds for every pair $e{\in}\s S^2M$.\eod
% [pd-rel hyperbolicity] @tag_4DFCC578:
\subsection{Relatively hyperbolic uniformities and the map theorem}\label{pd2floyd}
We propose a definition of relative hyperbolicity equivalent to the other known definitions
intended to the proof of the Floyd map theorem. This is a structure including a
``geometric part'' (graph) and ``dynamical part'' (divider and uniformity).
\deff{
$\s{RH_{pd}}$.
A uniformity $\Cal U$ on a connected $G$-set $M$ is called \it relatively hyperbolic \rm(we also
say that $\Cal U$ is a \it relative hyperbolicity \rm on $M$)
if it is generated by a perspective divider (see the definitions 3.2.3 and \ref{cmpMetrics})
and \it not parabolic \rm i.e, the $\Cal U$-boundary $\partial_{\Cal U}{=}\overline M{\setminus}M$
has at least two points.}

The metrics $\d_{v,\lambda}$, $v{\in}M$ defined in \ref{cmpMetrics}
determine on $M$ the same uniformity $\s U_{\Gamma,\lambda}$.
We call it the \it exponential Floyd uniformity\rm.
It is $G$-invariant.

 If $\lambda_1{\leqslant}\lambda_2$
then $\s U_{\Gamma,\lambda_1}{\subset}\s U_{\Gamma,\lambda_2}$.
If $\Gamma_1,\Gamma_2$ are different connecting structures for a connected $G$-set $M$
then the identity map is Lipschitz. Hence for every $\lambda$ one has
$\s U_{\Gamma_1,\lambda}{\subset}\s U_{\Gamma_2,\lambda^{1/C}}$ where $C$ is the maximum
of the $\Gamma_2$-length of the $\Gamma_1$-edges.
\setcounter{prop}5%
\begin{prop}[Map theorem]\label{mapThrm}
Let $G$ be a group, $M$ a connected $G$-set, $\Gamma$ a connecting graph structure
\hbox{\rm(see \ref{iDivider})} for $M$,
$\Cal U$ a relatively hyperbolic uniformity on $M$.
Then there exists $\lambda\in(0,1)$ such that $\Cal U$
is contained in the Floyd uniformity $\s U_{\Gamma,\lambda}$.
The inclusion induces a uniformly continuous
$G$-equivariant surjective map {\rm(called the \it Floyd map\rm)}
$(\overline M,\overline{\s U}_{\Gamma,\lambda})\to(\overline M,\overline{\Cal U})$
between the completions.\end{prop}\sl Proof\rm.
Since the uniformities $\s U_{\Gamma,\lambda}$ are $G$-invariant it suffices to prove that
every perspective divider $\b u$ belongs to some $\s U_{\Gamma,\lambda}$.

We fix $v{\in}M$ and define $B_n$ as in \ref{cmpMetr}.

Let $F$ be a finite set from \ref{idemp3}, let $\rho$ be such that $F\{v\}{\subset}B_\rho$ and let
$\lambda{\leftrightharpoons}2^{-1/\rho}$ as in \ref{floydConstants}.
By Frink lemma \ref{frinkLemma}, $\d_{\b u,F}^{-1}\left[0,{1\over4}\right){\subset}\b u$.
By \ref{cmpMetr},
there exists a constant $C{>}0$ such that $\d_{v,\lambda}^{-1}\left[0,{1\over4C}\right){\subset}\d_{\b u,F}^{-1}\left[0,{1\over4}\right)$.
Hence $\b u{\in}\s U_{\Gamma,\lambda}$.\eod
% [Karlsson] @tag_4DFA27C2:
\subsection{Generalized Karlsson lemma}\label{Karlsson}
The following plays a basic role in our theory of r.h. groups.
It is the main application of the map theorem \ref{mapThrm}.
\begin{prop}\label{KarlssonLemma}
Let $G$ be a group, $M$ a connected $G$-set, $\Gamma$ a connecting graph structure
\hbox{\rm(see \ref{iDivider})} for $M$, such that the action $G{\curvearrowright}\Gamma^1$
is proper \hbox{\rm(see \ref{properOnPairs})}.
Let $\Cal U$ be a relatively hyperbolic uniformity on $M$ \hbox{\rm(see \ref{pd2floyd})}.
For every entourage $\b v{\in}\Cal U$
there exists a finite set $E{\subset}\Gamma^1$ such that $\b v$
contains the boundary pair $\partial I$ of every geodesic segment $I$ that misses $E$.
\end{prop}Remarks.

1. We do not assume that the graph $\Gamma$ is locally finite.
In the proof we use the original Karlsson lemma \cite[Section 3, Lemma 1]{Ka03}
claiming that
 the Floyd length of a $\s d$-geodesic $I$ tends to zero while $\s d(v,I)\to\infty$.
It does not use local connectedness of the graph.

2. The following proof remains true after replacing the word `geodesic' by `quasigeodesic'
or even `$\alpha$-quasigeodesic' for a polynomial ``distortion function'' $\alpha$,
see \cite{GP11}. This is not used in the theory developed in this article
but is essential for many applications.
The discussion about Bowditch completion in \ref{BoCmpl} contains an example
of such application.
\vskip3pt
\sl Proof\rm.
Let $\b u$ be a perspective divider such that $\Cal U{=}\Cal U_{\b u}$.
By \ref{dividerEnt} we can assume that $\b v{=}{\cap}(S\{\b u\})$ for a finite $S{\subset}G$.
Let $F$ be a finite subset of $G$ containing $S$ and satisfying \ref{idemp3}.
Consider the dividing metric $\d_{\b u,F}$ from \ref{cmpMetrics}.
According to the notation of \ref{expansiveMetric} we have $\b{v{=}u}_2$.
Let $\lambda$ be a number from \ref{cmpMetr} and let $v{\in}M$ be a fixed ``reference'' vertex.

 By \ref{cmpMetr} and \ref{frinkLemma} there exists $\varepsilon_2{>}0$ such that $\d_{v,\lambda}\beta\geqslant\varepsilon_2$
for any $\beta{\notin}\b u_2$.
We also need $\varepsilon_3{>}0$ such that $\d_{v,\lambda}\beta\geqslant\varepsilon_3$
for any $\beta{\notin}\b u_3$.

By the original Karlsson lemma \cite{Ka03}
applied to $\d_{v,\lambda}$ there exists a number $d$
such that if $\s d(v,I){\geqslant}d$ then $\s{length}_{v,\lambda}I{<}\varepsilon_2$ and
thus $\partial I{\in}\b v$. Let $r$ be a number such that $\sum_{n\geqslant r-d}\lambda^n\leqslant\varepsilon_3$,
let $\b w{=}\root{2r}\of{\b u_3}$ and
let $E{\leftrightharpoons}\Gamma^1{\setminus}\b w$ be the set of all ``$\b w$-big'' edges.
By perspectivity of $\b w$ the set $E$ is finite.
We will verify that it satisfies the assertion.

Let $I$ be a geodesic segment whose edges do not belong to $E$.
Let $p$ be its vertex closest to the reference vertex $v$.
If $h{\leftrightharpoons}\s d(v,I){\geqslant}d$ then $\partial I{\in}\b v$ by the
choice of $d$. Suppose that $h{<}d$.
Let $J$ be the $r$-neighborhood of $p$ in $I$. So $I$ is the concatenation of three segments $I_-,J,I_+$.
We estimate the $\d_{v,\lambda}$-length of each $I_\pm{\leftrightharpoons}L$ provided that it is non-zero.
In this case $\s d(p,L){=}r$.
Let $x_r,x_{r+1},\dots$ denote the consecutive vertices of $L$ such that
$\s d(p,x_n){=}n$ for $n{\geqslant}r$.
By $\triangle$-inequality,
$$\s d(x_n,v)\geqslant\s d(x_n,p)-\s d(p,v)=
n{-}h\geqslant n-d\geqslant r-d.$$
Hence $\d_{v,\lambda}(x_n,x_{n+1})\leqslant\lambda^{n{-}d}$ and
$\s{length}_{v,\lambda}L{\leqslant}\varepsilon_3$, $\partial L{\in}\b u_3$.

Since the edges of $J$ belong to $\b w$ and $\s{length_d}J{\leqslant}2r$, we have $\partial J{\in}\b w^{2r}{\subset}\b u_3$.
Since $\b u_3^3{\subset}\b u_2$ we have $\partial I{\in}\b u_2{=}\b v$.\eod
% [floyd map] @tag_4A80F403) end of piece
% [visibility] @tag_4DFF2A59(
\section{Visibility}
% [Definition] @tag_4E006B73:
\subsection{Definition and examples}\label{visibility}
We now turn a weaker version of \ref{KarlssonLemma} into a definition.
Let $\Gamma$ be a connected graph with $\Gamma^0{=}M$.
For an edge $e{\in}\Gamma^1$ define
\begin{equation}\label{vizGenerators}\b u_e\leftrightharpoons\{\{x,y\}{\subset}\Gamma^0:
\s d(x,e)+\s d(y,e)\geqslant\s d(x,y)\}\end{equation}
This set consists in all pairs such that no geodesic segment joining them passes through $e$.
Since we do not allow multiple edges, the only edge that does not belong to $\b u_e$ is the edge $e$.

The filter $\s{Vis}\Gamma$ generated by $\{\b u_e:e{\in}\Gamma^1\}$ is the \it visibility filter \rm on $M$.
\deff{ \ref{visibility}
A uniformity $\Cal U$ on $M$ is called a \it visibility \rm on $\Gamma$ if
it is contained in $\s{Vis}\Gamma$.}
It follows from the remark after \ref{vizGenerators} that
any visibility is perspective in the sense of \ref{perspectivity}.
Hence if it exact then, by \ref{exactIsFine}, the graph $\Gamma$ is fine.

The generalized Karlsson lemma \ref{KarlssonLemma}
implies that, for any perspective divider $\b u$ for a connected $G$-set $M$
and for any connecting structure $\Gamma$ for $M$ the uniformity $\Cal U_{\b u}$ is a visibility on $\Gamma$.
Note that \ref{KarlssonLemma} actually claims something stronger:
each pair of vertices that can ``partially'' see each other outside of a big finite set of edges
is ``small''.

If $\Gamma$ is a locally finite (i.e, every vertex is adjacent to finitely many edges)
$\delta$-hyperbolic graph then $\s{Vis}\Gamma$ is a uniformity and hence the maximal possible visibility.

Every connected graph has the maximal visibility; we call it \it initial \rm since
there is a morphism from its completion to the completion of each other visibility.

If $\Gamma$ is locally finite then the filter, generated by the collection
of the sets of the form $\b v_E{\leftrightharpoons}\{$the pairs that can be joined by a path not containing
the edges in $E\}$, $E{\in}\s{Sub}^{<\infty}\Gamma^1$, is a visibility on $\Gamma$.
The completion with respect to it is just the Freudenthal's ``ends completion'' $\s{Fr}\Gamma$
and the boundary is the compact totally disconnected ``space of ends''.
(This fact is not used in the paper and leaved as an exercise.
It follows directly from our definitions and the notion of the ends of graphs,
see, for example, \cite{St71} or \cite{DD89}.)

We will see that many convergence group actions (actually all known) induce a visibility on the
Cayley graph of the acting group.

For a finite set $E{\subset}\Gamma^1$ we denote $\b u_E{\leftrightharpoons}{\cap}\{\b u_e:e{\in}E\}$.
The sets of the form $\b u_E$ are called \it principal\rm. They form a base for the filter $\s{Vis}\Gamma$.
\begin{prop}\label{vizCompact} Every visibility $\Cal U$ on a connected graph $\Gamma$ is pre-compact.\end{prop}\sl Proof\rm.
Let $\b{u{=}\sqrt v}$ for $\b v{\in}\Cal U$.
Since $\b u{\in}\s{Vis}\Gamma$ it contains a principal set $\b u_E$.
Let $S{\leftrightharpoons}{\cup}E{=}\{$the vertices of the edges in $E\}$.
For $x{\in}S$ the set
$\{v{\in}\Gamma^0:\s d(v,x){=}\s d(v,S)\}$ is $\b u^2$-small
since the shortest path from $v$ to $x$
does not contain the edges from $E$.
So $\Gamma^0$ is a union of finitely many $\b v$-small sets.\eod
\subsection{Visibility actions have Dynkin property}
The following proposition gives a large class of convergence group actions.
We do not know examples of convergence actions outside this class. However
it seems to be difficult to prove that every convergence action is a visibility action (see section \ref{geometric}).
If follows from a recent result of Mj Mahan \cite{Mj10}
that every Kleinian action of a finitely generated group is.
% [visIsDynkin] @tag_4E006B40:
\begin{prop}\label{visIsDynkin}Let a group $G$ act on a connected graph properly on edges and let $\Cal U$ be a $G$-invariant
visibility on $\Gamma$. Then the action has Dynkin property.
Furthermore the action on the completion has the convergence property.\end{prop}\sl Proof\rm.
Consider such an action $G{\on}\Gamma$.
By \ref{perspectivity} it is perspective,
so, for each entourage $\b u{\in}\Cal U$ and each finite set $S{\subset}\Gamma^0$
the set $\{g{\in}G:gS$ is not $\b u$-small\} is finite.

Let $\b{u,v}{\in}\Cal U$.
Since $\Cal U$ is a uniformity there exist principal sets $\b u_E$, $\b u_F$
such that $\b u_E^3{\subset}\b u$, $\b u_F^3{\subset}\b v$.
By the above remark it suffices to show that if
${\cup}E\in\s{Small}(\b u_F)$ and ${\cup}F\in\s{Small}(\b u_E)$
then $\b u{\bowtie}\b v$.

Let $A{\leftrightharpoons}\{x{\in}\Gamma^0:\s d(x,E){\geqslant}\s d(x,F)\}$, 
$B{\leftrightharpoons}\{x{\in}\Gamma^0:\s d(x,F){\geqslant}\s d(x,E)\}$.
Since a geodesic segment realizing the distance $\s d(x,F)$ ($x{\in}A$) can not contain
the edges from $E$ it is $\b u_E$-small. Thus $A$ is $\b u_E^3$-small hence $\b u$-small.
By the same reason $B$ is $\b v$-small.

The Dynkin property is verified.

By \ref{vizCompact} the $\Cal U$-completion $\overline{\Gamma^0}$ is compact,
by \ref{complDynkin} the action $G{\on}\overline{\Gamma^0}$
has Dynkin property, by \cite[5.3. Proposition P]{Ge09} it is discontinuous on triples i.e, has convergence
property (see {\ref{charProp}).
\eod
\vskip3pt
As a corollary we have the implication $\s{RH_{pd}{\Rightarrow}RH_{32}}$:
\begin{prop}\label{pd2fh}Let $\Cal U$ be a relative hyperbolicity on
a connected $G$-set $M$. Then the action on the completion $\overline M$
with respect to $\Cal U$ is relatively hyperbolic in the sense of $\s{RH_{32}}$.\end{prop}\sl Proof\rm.
By \ref{Karlsson} the uniformity $\Cal U$ is a visibility, by \ref{visIsDynkin} the action $G{\on}\overline M$
is discontinuous on triples.
It suffices to find a bounded ``fundamental domain'' for the action $G{\on}\b\Theta^2\overline M$.

Let $\b u$ be a perspective divider generating $\Cal U$ and let $\b v{=}\root3\of{\b u}$.
Let $\{\g{p,q}\}{\in}\b\Theta^2\overline M$.
By exactness of $\overline{\Cal U}$ (see \ref{Samuel}) there exists $\b w{\in}\Cal U$ such that
$\g{p{\cap}q}$
has no $\b w^3$-small sets.
For every pair $P{\in}\g p$, $Q{\in}\g q$ of $\b w$-small sets
we have $P{\times}Q\cap\b w=\varnothing$.

Since $\b u$ generates $\Cal U$ as a $G$-filter, $\b w$ contain a set ${\cap}(S\{\b u\})$
for $S{\in}\s{Sub}^{<\infty}G$.
Let $\b w_1{\leftrightharpoons}{\cap}(S\{\b v\})$ and let $P,Q$ be $\b w_1$-small sets in
$\g{p,q}$ respectively. Let $p{\in}P,q{\in}Q$.
Since $\{p,q\}{\notin}\b w$ there exists $g{\in}S$ such that $\{p,q\}{\notin}g\b u$.
We claim that the filter $g^{-1}\g p{\cap}g^{-1}\g q$ does not contain $\b v$-small sets.
Indeed if not, and $R$ is such a set, it would intersect $\b v$-small sets $g^{-1}P$ and $g^{-1}Q$
and hence $\{g^{-1}p,g^{-1}q\}{\in}\b u$. This is impossible.

So the complement of the entourage $\overline{\b v}$ (see \ref{Samuel}) is a bounded fundamental domain
for $G{\on}\b\Theta^2\overline M$.\eod
\vskip3pt
A relative hyperbolicity $\Cal U$ is not necessarily exact thus the completion map
$\iota_\Cal U:M\to\overline M$ can be not injective.
But $\Cal U$ is always a perspectivity for any connecting graph structure $\Gamma$
and hence $\Gamma$ can be replaced by another graph $\Delta$
as explained in \ref{perspectivityCompl}, on which $\Cal U$ induces an exact relative hyperbolicity.
We have the following.
\begin{prop}\label{exactening}Let $\Cal U$ be a relative hyperbolicity on a $G$-set $M$
and let $\Gamma$ be a connected graph with $\Gamma^0{=}M$ where $G$ acts properly on edges.
Then there exists a $G$-set $N$, an exact relative hyperbolicity $\Cal V$ on $N$,
a connected graph $\Delta$ with $\Delta^0{=}N$ and a uniformly continuous $G$-equivariant map
$\varphi:M\to N$ such that the induced map $\overline\varphi:\overline M\to\overline N$
is a homeomorphism.\qed\end{prop}
% [visibility] @tag_4DFF2A59)
% [alt-hyperbolicity] @tag_4DF2429A(
\section{Alternative hyperbolicity}\label{altHyp}
We are going to prove the equivalence between $\s{RH_{pd}}$ and the definition $\s{RH_{fh}}$ given in the introduction.
% [definition] @tag_4E010B81:
\subsection{Definition}
A connected graph $\Gamma$ is said to be \it alternatively hyperbolic \rm(we also say `alt-hyperbolic')
if the filter $\s{Vis}\Gamma$ is a uniformity and hence a visibility
(see \ref{visibility}).
This means that
\setcounter{equation}0%
\begin{equation}\label{altHyp1}\forall e\in\Gamma^1\
\exists
F\in\s{Sub}^{<\infty}\Gamma^1:\b u_F^2\subset\b u_e\end{equation}
where $\b u_F$ denotes the principal set derermined by $F$, see \ref{visibility}.
In turn this means that for every $e{\in}\Gamma^1$ there is finite set $F(e){\subset}\Gamma^1$
such that every geodesic triangle with $e$ on a side
contains an edge from $F(e)$ on another side.
\vskip3pt
Every locally finite $\delta$-hyperbolic graph is alt-hyperbolic with $F(e){=}\{f{\in}\Gamma^1:\s d(e,f){\leqslant}\delta\}$.
On the other hand if $\Gamma$ is alt-hyperbolic and $\{\s d(e,f):f{\in}F(e)\}$
is uniformly bounded then it is hyperbolic.

The classes \{hyperbolic graphs\} and \{alt-hyperbolic graphs\} are not included one to the other.
We will prove in \ref{fineHypIsAltHyp} that the hyperbolic fine graphs are alt-hyperbolic.
% [altRelHyp] @tag_4E01230F:
\subsection{Alternative relative hyperbolicity}\label{altRelHyp}
We make one more step towards the equivalence of $\s{RH_{pd}}$ to the other $\s{RH}$'s.
\deff{
$\s{RH_{ah}}$. An action of a group $G$ on a connected graph $\Gamma$ is
\it relatively hyperbolic \rm if $\Gamma$ is alt-hyperbolic and the action $G{\on}\Gamma^1$
is proper and cofinite and non-parabolic in the sense that no vertex is fixed by the whole $G$.}

To interpret $\s{RH_{ah}}$-action as an $\s{RH_{pd}}$ we only need to indicate a divider generating
the uniformity $\s{Vis}\Gamma$.
Let $E$ be a finite set of edges intersecting each $G$-orbit.
Then the entourage $\b u_E$ is a perspective divider.
The corresponding finite subset of $G$ from \ref{expansiveMetric}
is ${\cup}\{F(e):e{\in}E\}$ where $F(e)$ is given by \ref{altHyp1}.
The verification is straightforward.
% [fine hyp is alt] @tag_4DF09763:
\subsection{Hyperbolic fine graphs}
B. Bowditch noted that the Farb's ``conned-off'' graph for a relatively hyperbolic
group
(in the ``$\s{BCP}$ sense'' of Farb)
with respect to a collection $\Cal P$ of subgroups, is fine.
This was an important step in understanding the relative hyperbolicity.

\begin{prop}\label{fineHypIsAltHyp}
Every connected $\delta$-hyperbolic fine graph $\Gamma$
is alt-hyperbolic.
\end{prop}\sl Proof\rm.
This is an exercise on the common ``thin triangle'' techniques.
Inside the proof we locally change the style and the notation,
following \cite[Section 2.3]{GhH90}.

We regard graphs as $\s{CW}$-complexes, since we need points on the edges
other than the endpoints. Actually one additional point in each edge would suffice.

Let $M$ be a geodesic metric space.
We denote $|ab|{\leftrightharpoons}\s d_M(a,b)$.
The word `segment' will mean `geodesic segment'.
By $[ab]$ we denote a particular segment joining $a$ and $b$.
The reader should understand which of possible such segments we mean.

For $a,b,c\in M$ a \it triangle \rm$T$ is a union of
segments $S_a,S_b,S_c$ called \it sides \rm such that
$\partial S_a{=}\{b,c\}$, $\partial S_b{=}\{c,a\}$, $\partial S_c{=}\{a,b\}$.
The points $a,b,c$ are the \it vertices \rm of the triangle.
The sides are allowed to have non-trivial intersection.

A \it tripod \rm is a metric cone over a triple.
The points of the triple are the \it ends \rm and the ``vertex'' of the cone is its \it center\rm.
We regard a tripod as a triangle whose vertices are the ends.
For every triangle $T{=}abc$ there exists a unique
(up to isometry) \it comparison tripod \rm $T'{=}a'b'c'$ with center $t'$
and a \it comparison map \rm $x\mapsto x'$ taking isometrically the $T$-sides onto $T'$-sides.
It is \it short \rm, i.e, $|x'y'|{\leqslant}|xy|$ for all $x,y{\in}T$.

A triangle $T$ is $\delta$\it-thin \rm if
$|xy|{-}|x'y'|\leqslant\delta$ on $T$.
\vskip3pt
Now let $M$ be a graph $\Gamma$
and let $T{=}abc$ be a $\delta$-thin triangle in $\Gamma$ for a positive integer $\delta$.
Let $e{=}a_0b_0$ be an edge on the side $S_c$ not contained in $S_a{\cup}S_b$
such that $|ab|{=}|aa_0|{+}|a_0b_0|{+}|b_0b|$.
We will find a circuit $H$ in $\Gamma$ of length ${\leqslant}20\delta{+}6$
containing $e$ and having an edge in $S_a{\cup}S_b$.

Initially we construct the pieces of $H$ joining $a_0$ and $b_0$
with $S_a{\cup}S_b$.

Let $a_1$ be a vertex of $S_a{\cap}S_c$ closest to $a_0$.

In the \it exceptional \rm case `$|a_0a_1|{<}\delta$' denote $a_2{\leftrightharpoons}a_1$.
Otherwise
$a_2{\leftrightharpoons}$the vertex on $[a_1a_0]$ with $|a_2a_0|{=}\delta$.

Let $a_3$ be a vertex in $S_a{\cup}S_b$ closest to $a_2$. In the exceptional case
$a_2{=}a_3$. Otherwise $|a_2a_3|{\leqslant}\delta$ by $\delta$-thinness.
We choose and fix a segment $[a_2a_3]$.

Let $a_4$ be the vertex
of $[a_2a_3]{\cap}[a_2a_0]$ closest to $a_0$ (in the exceptional case $a_4{=}a_3{=}a_2{=}a_1$).
The arc $L_a{\leftrightharpoons}[a_0a_4]{\cup}[a_4a_3]$ (the thick line on the picture below)
of length ${\leqslant}2\delta$ joins
$a_0$ with $a_3$. By construction $L_a{\cap}(S_a{\cup}S_b){=}\{a_3\}$.

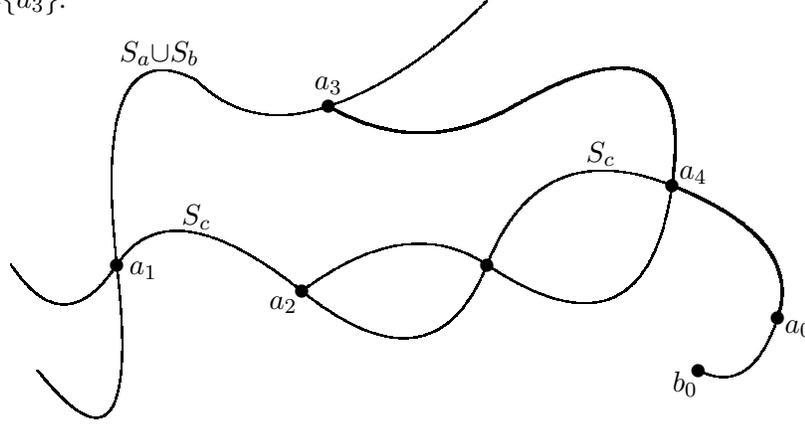
\begin{center}\begin{picture}(200,140)(-100,-53)
\renewcommand{\qbeziermax}{1000}
\put(-70,-10){\circle*{5}}% a_2
\put(-77,-15){\makebox(0,0)[c]{$a_2$}}
\qbezier(-70,-10)(-120,30)(-140,0)% S_c
\qbezier(-70,-10)(-20,-50)(0,0)
\qbezier(-70,-10)(-30,20)(0,0)
\put(0,0){\circle*{5}}
\qbezier(0,0)(60,-40)(70,30)
\qbezier(-110,70)(-90,50)(-60,60)
\qbezier(-110,70)(-150,90)(-140,0)
\put(-140,0){\circle*{5}}
\put(-130,-2){\makebox(0,0)[c]{$a_1$}}
\qbezier(-140,0)(-130,-90)(-170,-40)
\qbezier(-140,0)(-160,-30)(-180,0)
\qbezier(-60,60)(-30,70)(0,100)
\put(-124,80){\makebox(0,0)[c]{$S_a{\cup}S_b$}}
\put(-60,60){\circle*{5}}
\put(-60,68){\makebox(0,0)[c]{$a_3$}}
\qbezier(0,0)(20,50)(70,30)
\put(70,30){\circle*{5}}% a_4
\put(78,34){\makebox(0,0)[c]{$a_4$}}
\thicklines
\qbezier(70,30)(80,100)(10,60)
\qbezier(-60,60)(-25,40)(10,60)
\qbezier(70,30)(120,10)(110,-20)
\thinlines
\put(110,-20){\circle*{5}}% a_0
\put(118,-24){\makebox(0,0)[c]{$a_0$}}
\qbezier(110,-20)(100,-50)(80,-40)
\put(80,-40){\circle*{5}}% b_0
\put(75,-45){\makebox(0,0)[c]{$b_0$}}
\put(-110,18){\makebox(0,0)[c]{$S_c$}}
\put(43,42){\makebox(0,0)[c]{$S_c$}}
\end{picture}\end{center}

In the exceptional case $L_a{=}[a_0,a_3]{\subset}S_c$.

Replacing `$a$' by `$b$' we define the points $b_\iota$
for $\iota{\in}\overline{1,4}$ and the arc $L_b$.

Claim: $L_a{\cap}L_b{=}\varnothing$.
Indeed if both $L_a$ and $L_b$, are exceptional then they are disjoint segments
of $S_c$. If one of them is exceptional, say $L_b$, and the other is not
then both pieces $[a_0a_4]$ and $[a_4a_3]$ of $L_a$ are contained
in the $\delta$-neighborhood of $a_2$ which is disjoint from $[b_0b]$.
If both $L_a,L_b$ are not exceptional then they are contained in the disjoint
$\delta$-neighborhoods of $a_2$ and $b_2$ respectively.
So we have $L_a{\cap}L_b{=}\varnothing$ in all cases.

The set $L{\leftrightharpoons}L_a{\cup}e{\cup}L_b$ is an arc
of length ${\leqslant}4\delta{+}1$ with $L{\cap}(S_a{\cup}S_b){=}\partial L$.
The claim implies that $1{\leqslant}|a_3b_3|{\leqslant}4\delta{+}1$.

If $a_3$ and $b_3$ belong to $S{\in}\{S_a,S_b\}$ then the subsegment $[a_3b_3]$
of $S$ contains at least one edge,
has length ${\leqslant}4\delta{+}1$ and completes $L$ up to a circuit $H$
of length ${\leqslant}8\delta{+}2$.

Suppose now that $a_3$ and $b_3$ belong to different sides of $T$
and do not belong to $S_a{\cap}S_b$.

Let $t_a,t_b,t_c$ denote the preimages of the center $t'$ of $T'$
in the corresponding sides.
For $x{\in}S_c$ we have
\begin{equation}\label{altDistCenter}
|xt_c|\leqslant\s{max\{d}(x,S_a),\s d(x,S_b)\}
\end{equation}
Indeed if $x{\in}[at_c]$ then, for $p{\in}S_a$, the thinness inequality yields:
$|xt_c|{=}|x't'|{\leqslant}|x'p'|{\leqslant}|xp|$.

By \ref{altDistCenter}, the distance from $t_c$ to each of
$a_0,b_0$ is at most $2\delta{+}1$.
By the same reason the distance from each of $a_3,b_3$ to the
corresponding $t_?$ (which is $t_a$ or $t_b$) is at most $4\delta{+}1$.

Let $[a_3c]$ be the subsegment of the side in $\{S_a,S_b\}$ that contains $a_3$.
Similarly define $[b_3c]$.

Let $a_5{\in}[a_3c]{\cap}[b_3c]$ be the closest to $a_3$.
Since $a_5{\ne}a_3$ the segment $[a_5a_3]$ contains an edge of $S_a{\cup}S_b$.

If $|a_3a_5|{\leqslant}5\delta{+}2$ then
$H{\leftrightharpoons}L{\cup}[a_3a_5]{\cup}[a_5b_3]$
is a circuit of length ${\leqslant}2(5\delta{+}2{+}4\delta{+}1){=}18\delta{+}6$
with the desired properties.

So we can assume that $|a_3a_5|{>}5\delta{+}2$.

We actually repeat the construction of $L_a$ in the non-exceptional case: let
$a_6{\in}[a_3,a_5]$ be such that $|a_3a_6|{=}5\delta{+}2$.
Since $|a_3a_6|{\geqslant}|a_3't'|$ the comparison image $a_6'$ belongs to $[t'c']$
and $|a_6't'|{\geqslant}5\delta{+}2{-}(4\delta{+}1){=}\delta{+}1$.
Hence, by \ref{altDistCenter}, $\s d(a_6,S_c){\geqslant}|a_6't'|{\geqslant}\delta{+}1$.

Let $a_7$ be a point in the side containing $b_3$ closest to $a_6$.
By thinness, $|a_6a_7|{\leqslant}\delta$.

We choose a segment $I{=}[a_6,a_7]$.
The distance from each $x{\in}I$ to $S_c$ is ${\geqslant}1$
so it does not contain $e$.

Let $a_8{\in}[a_3a_6]{\cap}[a_6a_7]$ be the closest to $a_3$.
Since $|a_3a_8|{\geqslant}4\delta{+}2$ the segment $[a_3a_8]$
contains an edge in $S_a{\cup}S_b$.
The length of the circuit $H{\leftrightharpoons}L{\cup}[a_3a_8]{\cup}[a_8a_7]{\cup}[a_7b_3]$
is ${\leqslant}2(4\delta{+}1{+}5\delta{+}2{+}\delta){=}20\delta{+}6$.

So $\Gamma$ is alt-hyperbolic for $F(e){\leftrightharpoons}\{e\}{\cup}\{$the
edges of the circuits of length ${\leqslant}{20\delta{+}6}$ containing $e\}$.\qed
\vskip3pt
The proposition just proved shows that $\s{RH_{fh}{\Rightarrow}RH_{ah}}$.
% [alt-hyperbolicity] @tag_4DF2429A)
\section{More preliminaries}
% [general topology] @tag_4DFA802F(
The main purpose of the rest of the paper is the implication $\s{RH_{32}{\Rightarrow}RH_{pd}}$.
We give a proof under the following restriction: the $G$-set $M$ of non-conical points
is connected. This is trivially true for finitely generated groups.
In \cite{GP10} we will prove that $M$ is always connected.

There is also a ``non-main'' purpose:
to derive a general theory of convergence group actions preparing
some tools for the other theorems.
So our exposition is not ``absolutely minimal''.
We are trying to facilitate reading for those readers who are interested
only in our main purpose.

This section is a continuation of the preliminary Section \ref{preliminaries}.
Most of the information therein is widely known and can be found in the common sources.
% [actions] @tag_4A8351FE:
\subsection{Actions and representations}
Let $A,X,Y$ be sets.
Denote by $\s{Mp}(X,Y)$ the set of maps $X\to Y$. We consider an arbitrary map $\rho:A\to\s{Mp}(X,Y)$
as a \it representation \rm of a set $A$ by maps $X\to Y$ and as \it families \rm of maps $X\to Y$ indexed by $A$.
Our families of maps will be families of homeomorphisms. However this assumption is not always necessary.
In most cases  $A$ is a group, $X{=}Y$
(we refer to such cases as \it symmetric\rm), and $\rho$ is a homomorphism,
but sometimes we need to consider non-symmetric cases.

Denote by $\s{Bj}(X,Y)$ the set of bijective maps $X\to Y$.
For a family $\rho:A\to\s{Bj}(X,Y)$ a subproduct $U{\times}V$ of $X{\times}Y$ is \it invariant \rm if
$\rho(a)$ maps $U$ onto $V$ for each $a{\in}A$.

We can regard the maps $\alpha:A{\times}X\to Y$ as \it actions \rm of a set $A$ of ``operators''
on a set $X$ of ``points'' with values in another set $Y$.

Every set $K{\subset}X$ such that $\alpha(A{\times}K){=}Y$ is called \it generating\rm.
A generating set is sometimes called `fundamental domain' for the action.

The ``exponential law'' $\s{Mp}(A,\s{Mp}(X,Y))\simeq\s{Mp}(A{\times}X,Y)$
is a natural bijection\begin{center}\{representations$\}{\leftrightarrow}\{$actions\}.\end{center}
If $\rho$ is a representation and $\alpha$ is an action then by $\rho_*,\alpha_*$ we denote
respectively the corresponding action and the corresponding representation.
In the symmetric case $\alpha$ is a \it group action \rm if and only if $\alpha_*$ is a homomorphism of groups.

In case when our sets $A,X,Y$ are topological spaces we suppose by default that the actions are \bf continuous\rm.
However sometimes we have to prove the continuity of an action given by some construction.
When an action $\alpha$ is continuous then $\alpha_*A$ is contained in the set $\s{Top}(X,Y)$ of continuous maps $X\to Y$.
Moreover the representation $\alpha_*$ is continuous with respect to certain topology on $\s{Top}(X,Y)$.
In case when $X$ is locally compact (which is always the case in this paper) this topology is compact-open.

Our default topology on $\s{Top}(X,Y)$ is \it compact-open\rm.
It is well-known (see, e.g. \cite[Theorem \romannumeral12.3.1]{Du66}, \cite[\romannumeral7.8]{McL98})
that if $X$ is locally compact then the exponential law gives
one-to-one correspondence between continuous actions and continuous representations:
\begin{equation}
\s{Top}(A,\s{Top}(X,Y))\simeq\s{Top}(A{\times}X,Y)
\end{equation}

For every set $A{\subset}\s{Top}(X,Y)$ the inclusion map can be regarded as a representation.
The corresponding action of the space $A$ is an \it evaluation \rm action.

An action $\alpha:A{\times}X\to Y$ is \it cocompact \rm if
it possesses a compact generating set.
% [morphism] @tag_4E031BBB:
\subsection{Morphisms}\label{morphism}
Let $\rho_\iota:L_\iota\to\s{Homeo}(X_\iota,Y_\iota)$ for $\iota\in\{0,1\}$ be continuous locally compact families
of homeomorphisms between compactums.
A \it morphism \rm$\rho_0\to\rho_1$ is a triple
\setcounter{equation}0
\begin{equation}\label{morph}
(L_0\overset\alpha\to L_1,X_0\overset\beta\to X_1,Y_0\overset\gamma\to Y_1)\end{equation}
of continuous maps with $\alpha$ \bf proper \rm
respecting the actions in a natural way.
If $\mu:A{\times}X\to Y$ is an action then $(\mu_*,\s{id}_X,\s{id}_Y)$ is the \it tautological \rm morphism from $\mu$
to the evaluation action $(\mu_*A){\times}X\to Y$.
\subsection{Group actions}
For a topological group $G$ a $G$\it-space \rm is a topological space $T$ where $G$
acts continuously by homeomorphisms.
A subspace of a $G$-space $T$ is is $G$\it-bounded \rm if its image in the quotient
space $T/G$ is bounded. In the discrete case $G$-bounded sets are $G$\it-finite\rm.
% [compactifying] @tag_4DF335C4:
\subsection{Compactification of a locally compact family}\label{compactification}
For a locally compact space $B$ denote by $\widehat B$ the one-point compactification $B{\cup}\{\infty_B\}$.

For a closed subset $R$ of a compactum $S$ denote by $S/R$ the result of collapsing $R$ to a point.
On the categorical language
$S/R$ is the pushout space of the diagram $1\leftarrow R\hookrightarrow S$.
The immediate consequence of the above consideration is the following
\begin{prop}\label{remainder}
Let $S$ be a compactum, $R{\in}\s{Closed}S$ and $f:B\to S{\setminus}R$ be a continuous map
from a locally compact Hausdorff space $B$.
The following properties are equivalent:\hfil\penalty-10000
$\s a:$ $f$ is proper;\hfil\penalty-10000
$\s b:$ the map $\widehat f:\widehat B\to S/R$ that maps $\infty_B$ to the point $R$
is continuous.
\end{prop}
A continuous map $f:B\to S$ is called $R$\it-compactifiable \rm if $\varnothing{=}R{\cap}fB$ and,
as a map $f:B\to S{\setminus}R$, it satisfies conditions $\s a$--$\s b$ of \ref{remainder}.
For such a map the top horizontal arrow of the pullback square
\setcounter{equation}1
\vskip10pt
\begin{equation}\label{rComp}%\hbox{\vbox to 40pt{\vfill}}
\vcenter{\begin{picture}(50,30)(-220,-30)
%\put(0,21){\makebox(0,0)[b]{$\scriptstyle f+R$}}
\put(-20, 20){\makebox(0,0){$C$}}
\put( 20, 20){\makebox(0,0){$S$}}
\put(-14,20){\vector(1, 0){27}}
\put(-20,15){\vector(0,-1){28}}
\put( 20,15){\vector(0,-1){28}}
%\put(-21, 0){\makebox(0,0)[r]{$\scriptstyle g$}}
\put( 21, 0){\makebox(0,0)[l]{$\scriptstyle\s{pr}$}}
\put(- 5,-17){\makebox(0,0)[b]{$\scriptstyle\widehat f$}}
\put(-20,-20){\makebox(0,0){$\widehat B$}}
\put(-15,-20){\vector(1, 0){23}}
\put( 20,-20){\makebox(0,0){$S/R$}}
\end{picture}}
\end{equation}
 is called the $R$\it-compactification\rm. Denote it by $f{+}R$.
The corresponding space $C$ is the union of homeomorphic copies of $B$ and $R$.
Thus we denote it by $B{+}_fR$.
\vskip3pt
The following is an easy consequence of the definitions.
\begin{prop}\label{RFcmp}
Let $S$, $R$, $B$ and $f$ be as in \ref{remainder}.
If $F{\in}\s{Closed}S$ and $\overline{fB}{\setminus}fB{\subset}F{\subset}R$
then $f$ is $R$-compactifiable if and only if it is $F$-compactifiable.\end{prop}

% [vietoris] @tag_4A9C6F62
\subsection{Vietoris topology}\label{Vietoris}
As in \cite{Ge09} we need certain facts about the Vietoris topology.
Most of them are simple exercises in general topology.

Let $Z$ be a compactum. On the set $\mathcal Z{\leftrightharpoons}\mathsf{Closed}Z$
consider the topology $\mathsf{OG}$ defined by declaring \bf open \rm the sets of the form
$o^\downarrow{\leftrightharpoons}\mathcal Z{\cap}\s{Sub}(o)$ ($o{\in}\mathsf{Open}Z$);
and the topology $\mathsf{CG}$ defined by declaring \bf closed \rm
the sets of the form $c^\downarrow$ ($c{\in}\mathcal Z$).
The sum $\mathsf{Vi\leftrightharpoons OG{+}CG}$ is the \it Vietoris topology\rm.
It will be our default topology on $\mathcal Z$.

The space $\mathcal Z$ is a compactum.
If $\b u{\in}\s{Ent}(Z)$ then $\mathsf{Vi}(\b u)\leftrightharpoons\{\{c,d\}{\in}\s S^2\mathcal Z:
c{\subset}d\b u,d{\subset}c\b u\}\in\s{Ent}\mathcal Z$.
The operator $\b u\mapsto\mathsf{Vi}(\b u)$ preserves inclusions
and maps cofinal subsets of the filter
$\mathsf{Ent}Z$ to cofinal subsets of the filter $\mathsf{Ent}\mathcal Z$.

The union map $\mathcal Z^2\to\mathcal Z$,
$(c_0,c_1)\mapsto c_0{\cup}c_1$ is continuous with respect to the topologies
$\mathsf{OG}$, $\mathsf{CG}$ and hence with respect to $\mathsf{Vi}$.

Every continuous map $f:Z_0\to Z_1$ between compactums induces the map $f_*:\s{Closed}Z_0\to\s{Closed}Z_1$ continuous with respect
to each of the three topologies.
If $f$ is injective then the map $f^*:\s{Closed}Z_0\leftarrow\s{Closed}Z_1$ of taking preimage
is $\s{OG}$-continuous, but not necessarily $\s{Vi}$-continuous.
\vskip3pt
Let $Z$ and $\mathcal Z$ be as above.
\vskip3pt
\bf Lemma\sl. The set
$\Delta{\leftrightharpoons}\{\tupl{p,q}{\in}\mathcal Z^2:p{\supset}q\}$
is closed in the topology $\s{OG{\times}CG}$ on $\mathcal Z^2$.\hfil\penalty-10000
\it Proof\rm.
If $p{\not\supset}q$ then $\exists\frak q{\in}q{\setminus}p$.
Let $o,v$ be disjoint open neighborhoods of $p$ and $\frak q$ respectively.
Then $o^\downarrow{\times}((v')^\downarrow)'$
is an $\s{OG{\times}CG}$-neighborhood of the point $\tupl{p,q}{\in}\mathcal Z^2$
disjoint from $\Delta$.\eod
\vskip3pt
For a set $\mathcal{C{\subset}Z}$ let $\mathcal C^{\{\downarrow\}}{\leftrightharpoons}{\cup}\{c^\downarrow:c{\in}\mathcal C\}$.
\vskip3pt
 \bf Corollary\sl.
For any $\mathcal{C{\in}\s{Closed}Z}$ the set $\s{Loc_{\tupl{\mathcal Z,OG}}\mathcal C{\cap}Closed_{CG}}\mathcal Z$
generates the filter $\s{Loc_{Vi}}(\mathcal C^{\{\downarrow\}})$.\hfil\penalty-10000
\it Proof\rm.
If $x{\notin}\mathcal C^{\{\downarrow\}}$ then $\varnothing=\Delta\cap\mathcal C{\times}\{x\}$.
Since $\Cal C$ is compact (with respect to $\s{Vi}$ and hence with respect to weaker topologies $\s{OG}$ and $\s{CG}$)
by Lemma and the Walles Theorem \ref{Walles}, there exist
a $\s{OG}$-open $\mathcal{A{\supset}C}$ and a $\s{CG}$-open $\mathcal B{\supset}\{x\}$
such that
$\varnothing=\Delta\cap\mathcal{A{\times}B}$.
Thus $\mathcal B'$ is a $\s{CG}$-closed $\s{OG}$-neighborhood of $\mathcal C$ and hence of $\mathcal C^{\{\downarrow\}}$
that does not contain $x$.
The result follows from \ref{closedGenFil}.\eod
\vskip3pt
For a collection $\mathcal S$ of sets denote by $\s{min}\mathcal S$
the set of all minimal elements of $\mathcal S$.
\begin{prop}\label{OG}
If $\mathcal{C,S}{\in}\s{Closed_{Vi}}\mathcal Z$, $\mathcal C{\subset}\s{min}\mathcal S$
then every $\s{Vi}$-neighborhood of $\mathcal C$ in $\mathcal S$
contains an $\s{OG}$-neighborhood of $\mathcal C$. That is
$\s{Loc}_{\tupl{\mathcal S,\s{Vi}}}\mathcal C=
\s{Loc}_{\tupl{\mathcal S,\s{OG}}}C$.\hfil\penalty-10000
\end{prop}
\it Proof\rm.
For a $\s{Vi}$-open neighborhood $\mathcal V$ of $\mathcal C$ in $\mathcal S$
the $\s{Vi}$-open subset $\mathcal V{\cup}\mathcal S'$ of $\mathcal Z$
contains $\mathcal C^{\{\downarrow\}}$ since $\mathcal{S{\cap}C^{\{\downarrow\}}{=}C}$.
By Corollary, $\mathcal V{\cup}\mathcal S'$ contains some
$\mathcal{W{\in}\s{Loc_{OG}}(C}^{\{\downarrow\}})$,
so we have $\mathcal{S{\cap}W\subset S{\cap}V}$ as required.
\eod
\vskip5pt
The group $\s{Homeo}(Z,Z)$ acts naturally on $\mathcal Z{=}\s{Closed}Z$. This action is
continuous since if a homeomorphism is uniformly $\b u$-close to the identity map then
$\{c,\varphi c\}{\in}\s{Vi}(\b u)$ for every $c{\in}\mathcal Z$. Taking into account the
correspondence between actions and representations we obtain
\begin{prop}\label{viCont} If a topological group acts continuously on a compactum $Z$ then
the induced action on the space $\s{Closed}Z$ is continuous.\end{prop}
% [qHomeo] @tag_4AA7C6BB:
\subsection{Graphical embeddings and quasihomeomorphisms}\label{qHomeo}
Let $X,Y$ be compactums and let $Z{\leftrightharpoons}X{\times}Y$, $\mathcal Z{\leftrightharpoons}\mathsf{Closed}Z$.
For $\tupl{\b{u,v}}\in\s{Ent}X\times\s{Ent}Y$ the set\hfil\penalty-10000
\centerline{$\b{u{\cdot}v}\leftrightharpoons\{\{\tupl{\frak p_0,\frak q_0},\tupl{\frak p_1,\frak q_1}\}:
\{\frak p_0,\frak p_1\}{\in}\b u,\{\frak q_0,\frak q_1\}{\in}\b v\}$}\hfil\penalty-10000
is an entourage of $Z$.
The operator $\b{\tupl{u,v}\mapsto u{\cdot}v}$ preserves inclusions.
If $\mathcal P$ is cofinal in $\s{Ent}X$ and $\mathcal Q$ is cofinal in $\s{Ent}Y$ then
$\{\b{u{\cdot}v:\tupl{u,v}}{\in}\mathcal{P{\times}Q}\}$ is cofinal in $\s{Ent}Z$.

The set $\mathcal S{\leftrightharpoons}\mathsf{Surj}Z$
of all surjective (see \ref{composition}) \bf closed \rm subsets is $\mathsf{OG}$-closed
(hence $\mathsf{Vi}$-closed) in $\mathcal Z$.

By assigning to a continuous map $f:X\to Y$ its \it graph \rm$\boldsymbol\gamma f{\leftrightharpoons}\{\tupl{\frak p,f\frak p}:\frak p{\in}X\}$
we obtain the \it graphical embedding \rm$\mathsf{Top}(X,Y)\to\mathcal Z$
which is a homeomorphism onto its image. We identify continuous maps with
its graphs. So the set $\mathcal H{\leftrightharpoons}\mathsf{Homeo}Z{\leftrightharpoons}\mathsf{Homeo}(X,Y)$
of homeomorphisms is a subset of $\s{min}\mathcal S$ (see \ref{Vietoris}).

The closure $\overline{\mathcal H}$ of $\mathcal H$ in $\mathcal Z$
is contained in $\mathcal S$ but not in $\s{min}\mathcal S$ in general.
The points of $\overline{\mathcal H}$
are the \it quasihomeomorphisms \rm from $X$ to $Y$.
For a set $H{\subset}\mathcal H$ the points of the \it remainder \rm$\overline H{\setminus}H$
are the \it limit quasihomeomorphisms \rm of $H$.

A subset $s$ of $Z$ is \it single-valued \rm at a point $\frak p{\in}X{\sqcup}Y$ if the intersection of $s$ with
the fiber over $\frak p$ is a single point.
A closed set $s$ belongs to $\mathcal H$ if and only if it is single-valued at every point of $X{\sqcup}Y$.

\begin{prop}\label{eqCont}
A set $\Phi{\subset}\mathcal H$ is equicontinuous at a point $\frak p{\in}X$ if and only if
every $s{\in}\overline\Phi$ is single-valued at $\frak p$.\end{prop}
\it Proof\rm.
Assume that $\Phi$ is equicontinuous at $\frak p$.
Let $\tupl{\frak{p,q}_0},\tupl{\frak{p,q}_1}\in s{\in}\overline\Phi$.

For an entourage $\b v{\in}\s{Ent}Y$ there exists a neighborhood $p$ of $\frak p$
such that $\varphi p$ is $\b v$-small for all $\varphi{\in}\Phi$.
Let $\b u$ be an entourage of $X$ such that $\frak p\b u{\subset}p$ and let
$\varphi{\in}\Phi$ be a homeomorphism contained in the $\b{u{\times}v}$-neighborhood of $s$ (see \ref{Vietoris}).
Let $\tupl{\frak p_0,\varphi\frak p_0},\tupl{\frak p_1,\varphi\frak p_1}$ be points $\b{u{\times}v}$-close to
$\tupl{\frak{p,q}_0}$ and $\tupl{\frak{p,q}_1}$ respectively.
We have $\{\frak q_0,\frak q_1\}{\in}\b v^3$. Since $\b v$ is an arbitrary entourage of $Y$ we have $\frak q_0{=}\frak q_1$.

Assume that $\Phi$ is not equicontinuous at $\frak p$. Thus there exists an entourage $\b v$ of $Y$ such that
for each neighborhood $p$ of $\frak p$ the set $\Phi_p{\leftrightharpoons}
\{\varphi{\in}\Phi:\varphi p\notin\s{Small}(\b v^3)\}$ is nonempty.
We have a filtered family $p\mapsto\Phi_p$ in a compactum $\mathcal Z$. So it possesses an accumulation point $s$.
The set $q{\leftrightharpoons}s(\frak p)$ is not $\b v$-small since $q\b v$ contains a non-$\b v^3$-small pair.\eod
\vskip3pt
Remark. The proof of \ref{eqCont}
remains valid in a more general situation when $\Phi$ is a set of continuous maps $X\to Y$
(not necessarily homeomorphisms).
We do not need this generalization.
\begin{prop}\label{evCont}Let $S$ be a closed set of quasihomeomorphisms single-valued at a point $\frak p{\in}X$.
Then the natural evaluation map $\s{ev}_{\frak p}:S{\ni}s\mapsto s\frak p{\in}Y$ is continuous.\end{prop}
\it Proof\rm. If $f{\in}\s{Closed}Y$ then $\s{ev}_{\frak p}^{-1}f{=}\{s:s\cap\frak p{\times}f\ne\varnothing\}$ is $\s{OG}$-closed.\qed
% [topology] @tag_4DFA802F) end of range
% [x-property] @tag_4AA83185(
\section{Convergence property}
% [crosses] @tag_4AA831B0:
\subsection{Crosses and their neighborhoods}\label{cross}
In this section $X,Y,Z,\mathcal{Z,S,H}$ mean the same as in \ref{qHomeo}.

For $\frak{p{=}\tupl{r,a}}{\in}Z$ the set
$\frak{p^\times\leftrightharpoons r}{\times}Y\cup X{\times}\frak a=(\frak{r'{\times}a'})'$ is a \it cross \rm with
\it center \rm$\frak p$. The point $\frak r$ is the \it repeller \rm and $\frak a$ is the \it attractor \rm for $\frak p^\times$.
A cross $\tupl{\frak{r,a}}^\times$ is said to be \it diagonal \rm if $\frak{r{=}a}$.

The \it cross map \rm$\zeta_Z:\frak{p\mapsto p}^\times$ is a continuous map $Z\to\mathcal Z$ (see \ref{Vietoris}).
So its image $\s{Im}\zeta_Z$ is a closed subset of $\mathcal S$.
A cross is not a minimal element of $\mathcal S$ if and only if $|Z|{>}1$ and its center is an isolated point of $Z$.
So the set $\mathcal C{\leftrightharpoons}\{$the crosses with nonisolated center\} is closed in $\mathcal S$.
It is contained in $\s{min}\mathcal S$.
Note that it is empty if $Z$ is finite.

\begin{prop}\label{nonIsol}
 If $|Z|{>}1$ and a cross $c$ is a quasihomeomorphism then its repeller and attractor
are nonisolated in $X$ and $Y$ respectively. In particular $c$ belongs to $\mathcal C$.
If $|X|{\geqslant}3$ and a cross $c$ contains a quasihomeomorphism $q$ then $c{=}q$.\end{prop}
\it Proof\rm.
If $Z$ is finite then $\mathcal Z$ is also finite and the topology $\s{Vi}$ is discrete.
In particular any quasihomeomorphism is a homeomorphism.
Since $|Z|{>}1$ any cross is not a homeomorphism and hence not a quasihomeomorphism.
If $|X|{\geqslant}3$ then no homeomorphism is contained in a cross.

Suppose $|Z|{=}\infty$.

Let $c{=}\tupl{\frak{r,a}}^\times$. Assume that $\frak r$ is isolated.
Let $a$ be a neighborhood of $\frak a$ such that $a'{\leftrightharpoons}Y{\setminus}a$ is not a single point.
Then the neighborhood $\{\frak r\}{\times}Y\cup X{\times}a=(\frak r'{\times}a')'$ of $c$ and hence of $q$ does not contain a homeomorphism.

Thus $\frak{r,a}$ are not isolated hence $c$ is a minimal surjective correspondence. Hence $c{=}q$.\eod
\vskip5pt

A set $\Phi{\subset}\mathcal H$ is $\times$\it-compactifiable \rm if $\Phi{\cup}\mathcal C$ is closed.
In terms of \ref{compactification} `$\times$-compactifiable' means that the inclusion $\Phi\hookrightarrow\mathcal S$
is a $\mathcal C$-compactifiable family.

Since $\mathcal H{\subset}\s{min}\mathcal S$,
it follows from \ref{OG} that each neighborhood in $\mathcal S$ of any closed subset $\mathcal D$ of $\mathcal{H{\cup}C}$
contains an $\s{OG}$-neighborhood. In particular the topologies on $\mathcal D$ induced by $\s{Vi}$ and by $\s{OG}$ coincide.

Now \ref{Walles} implies

\begin{prop}\label{locCross}
Every neighborhood in $\mathcal S$ of a cross $\tupl{\frak{r,a}}^\times{\in}\s{min}\mathcal S$
contains a neighborhood of the form $((r'{\times}a')')^\downarrow$ where $r{\in}\s{Loc}_X\frak r$ and
$a{\in}\s{Loc}_Y\frak a$.
\end{prop}
Consider now the composition \ref{composition} of crosses as binary correspondences.
Note that\hfil\penalty-10000
$\frak{\tupl{a,b}^\times\circ\tupl{c,d}^\times=\tupl{a,d}^\times}$ if $\frak{b{\ne}c}$.
\begin{prop}\label{crossAlgebra}
The composition \ref{composition} of binary correspondences is continuous
at each point $p\leftrightharpoons\frak{\tupl{\tupl{a,b}^\times,\tupl{c,d}^\times}}$ such that $\frak{b{\ne}c}$
 and $\tupl{\frak{a,d}}$ is not isolated in $X{\times}U$.\end{prop}
\it Proof\rm.
In view of \ref{locCross} consider a neighborhood
$\mathcal U{\leftrightharpoons}((a'{\times}d')')^\downarrow$ where $a{\in}\s{Loc}_X\frak a$, $d{\in}\s{Loc}_U\frak d$.
If $b,c$ are disjoint neighborhoods of $\frak b$ and $\frak c$ respectively then
$$(a{\times}Y\cup X{\times}b)\circ(c{\times}U\cup Y{\times}d)=(a{\times}U\cup X{\times}d)$$ hence the composition map
takes the neighborhood $((a'{\times}b')')^\downarrow{\times}((c'{\times}d')')^\downarrow$ of $p$ to $\mathcal U$.\eod
% [char properties] @tag_4A9BCD18:
\subsection{Characteristic properties of convergence actions}\label{charProp}
Let $X,Y$ be compactums of cardinality ${\geqslant}3$.
The following properties of a locally compact family
$\rho:A\to\s{Homeo}(X,Y)$ of homeomorphisms
are equivalent \cite[Proposition P, subsection 5.3]{Ge09}:\hfil\penalty-10000
$\s a:$ $\rho$ is 3-proper, i.e. the induced family
$\mathbf\Theta^3\rho:A\to\s{Homeo}(\mathbf\Theta^3X,\mathbf\Theta^3Y)$ is proper;\hfil\penalty-10000
$\s b:$ $\rho$ is $\times$\it-compactifiable \rm i.e., $\mathcal C$-compactifiable in the sense of \ref{remainder}
as a map $A\to\mathcal Z$;\hfil\penalty-10000
$\s c:$ $\rho$ is \it Dynkin \rm i.e., for every entourages $\b u{\in}\s{Ent}X$, $\b v{\in}\s{Ent}Y$
the set $\{\frak a{\in}A:\b v{\#}\frak a^\rho\b u\}$ is bounded in $A$.

A family satisfying the properties $\s a$--$\s c$ is called a $\times$\it-family \rm or a $\times$\it-representation\rm.
The corresponding action is a $\times$\it-action \rm or an action with the \it convergence property\rm.

Remark 1. If $A$ is not compact then the set $\rho A$ possesses limit crosses. Thus $X$ contains
nonisolated points, i.e,  $|X|{=}\infty$.

Remark 2. If follows from \ref{RFcmp} that `$\mathcal C$-compactifiable' family of homeomorphisms is
the same as a family compactifiable by the set $Z^{\{\times\}}{\leftrightharpoons}\{$crosses$\}$.
We will consider the compactifications described in \ref{compactification}.
% [limit sets] @tag_4A9BCD2C:
\subsection{Limit set operators}\label{limitSet}
For a $\times$-action $\mu:L{\times}X\to Y$ and a set $S{\subset}L$ denote by $\partial_\mu S$
the set of the centers of the
limit crosses of a subset $S{\subset}L$ and by $\partial_{\mu,0}S$, $\partial_{\mu,1}S$ the projections of
$\partial_\mu S$ onto $X$ and $Y$ respectively.
We call $\partial_\mu S$, $\partial_{\mu,0}S$, $\partial_{\mu,1}S$ respectively the $\times$\it-remainder \rm of $S$,
the \it repelling set\rm, and the \it attracting set\rm.
We sometimes omit the index `$\mu$' writing $\partial,\partial_0,\partial_1$ for these operators.
\begin{prop}\label{limSet}
Let $\mu:L{\times}X\to Y$ be a $\times$-action and let $W{=}U{\times}V\in\mathcal Z$ be an invariant closed subproduct.
If\, $W{\supset}\partial_\mu L$ then the restriction $\lambda$ onto $U'{\times}V'$ is proper.
If $|U|{\geqslant}2$ then $\partial_\mu L{\subset}W$.
If $|U|{\geqslant}3$ then $\nu{\leftrightharpoons}\mu|_{L\times U}$ is an $\times$-action
and $\partial_\mu L{=}\partial_\nu L$.
\end{prop}
\it Proof\rm.
We can assume that $L$ is not compact thus the $\times$-remainder $\partial_\mu L$ is nonempty.
In particular there is at least one homeomorphism $U\to V$ and hence $|V|{=}|U|$.
Let $K{\subset}U',K_1{\subset}V'$ be compact subsets.
If $S{\leftrightharpoons}\{g{\in}L:gK{\cap}K_1{\ne}\varnothing\}$
is not bounded in $L$ then it possesses a limit cross $\tupl{\frak{r,a}}^\times$.
On the other hand, $X{\setminus}K{\in}\s{Loc}_X\frak r$ and $Y{\setminus}K_1{\in}\s{Loc}_Y\frak a$
thus there exists $g{\in}S$ such that $gK{\subset}Y{\setminus}K_1$. A contradiction.

Suppose that $|U|{\geqslant}2$.
Let $\tupl{\frak{r,a}}{\in}\partial_\mu L$ and let $\frak p{\in}U{\setminus}\frak r$.

For an arbitrary neighborhood $a$ of $\frak a$ there exists a homeomorphism $\varphi{\in}\mu_*L$ contained
in the neighborhood $(\frak p{\times}a')'$ of the cross $\tupl{\frak{r,a}}^\times$. So $\varphi\frak p{\in}V{\cap}a$.
So $\frak a{\in}\overline V{=}V$.

By the same reason we have $\frak r{\in}U$.

If $|U|{\geqslant}3$ then the restriction $\nu$ of a 3-proper action $\mu$ is obviously 3-proper.
The restriction map $\mathcal{Z\to W}{\leftrightharpoons}\s{Closed}W$ maps limit crosses for $\mu$ to crosses.
It is $\s{OG}$-continuous, see \ref{Vietoris}.
Thus the restriction over the compactum $\mu_*L\cup\mathcal C$ is continuous.
Hence it maps the closure surjectively onto the closure of the image.\eod
\vskip3pt
\bf Corollary\sl. Let $\mu:G{\times}T\to T$ be a $\times$-action of a locally compact group on a compactum $T$
and let $U$ be a closed invariant subset of $T$ containing at least two points.
Then the \it limit set \rm$\bold\Lambda G{\leftrightharpoons}\partial_0G{=}\partial_1G$
is contained in $U$\rm.
\vskip3pt
The points of the limit set are the \it limit points \rm of the action.
\vskip5pt
For a $\times$-action $\mu:L{\times}X\to Y$, an entourage $\b v{\in}\s{Ent}Y$ and $F{\in}\s{Closed}X$ put\hfil\penalty-10000
$\s{Big}_\mu(F,\b v){\leftrightharpoons}\{g{\in}L:gF{\not\subset}\s{Small}(\b v)\}$.

\begin{prop}\label{big} $\partial_0\s{Big}_\mu(F,\b v){\subset}F$\end{prop}
\it Proof\rm. Let $\tupl{\frak{r,a}}^\times$ be a limit cross for $\s{Big}_\mu(F,\b v)$.
Assume that $\frak r{\notin}F$. Let $a$ be a $\b v$-small neighborhood of $\frak a$.
For some $g{\in}\s{Big}_\mu(F,\b v)$ the homeomorphism $\mu_*g$ is contained in the neighborhood
$F'{\times}Y{\cup}X{\times}a$ of $\tupl{\frak{r,a}}^\times$.
It maps $F$ into a $\b v$-small set $a$. A contradiction.\eod
% [elementary] @tag_4AB7C2BD:
\subsection{Elementary $\times$-actions}\label{elementary}
Let $\mu:G{\times}T\to T$ be a $\times$-action of a locally compact group $G$.
It follows from \ref{crossAlgebra} that the $\times$-remainder $\partial G$ is closed under the \bf partial \rm operation
\setcounter{equation}0
\begin{equation}\label{crossOp}
\frak{\tupl{a,b}{\cdot}\tupl{c,d}{\leftrightharpoons}\left\{{\displaystyle\tupl{a,d}\text{ if }b{\ne}c\atop
\displaystyle\text{undefined if }b{=}c}\right.}\end{equation}
and \it symmetric \rm i.e. invariant under the transposition map $\frak{\tupl{p,q}\mapsto\tupl{q,p}}$.

The following ``algebraic'' lemma describes the sets of this types.
\setcounter{prop}1
\begin{prop}\label{xRemainder}
Let $T$ be a set and let $D$ be a symmetric subset of $T^2$ closed under the partial operation \ref{crossOp}.
Then one of the following is true:\hfil\penalty-10000
$\s a:$ $D$ has the form $\{\frak{\tupl{p,q},\tupl{q,p}}\}$;\hfil\penalty-10000
$\s b:$ there exists $\frak p{\in}T$ such that $\tupl{\frak{p,p}}{\in}D{\subset}\tupl{\frak{p,p}}^\times$;\hfil\penalty-10000
$\s c:$ $D{=}L^2$ for $L{\subset}T$.\end{prop}
\it Proof\rm. Suppose that $D$ is not contained in a diagonal cross. To prove that $D$ is of type $(\s c)$
it suffices to prove that $\tupl{\frak{a,b}}{\in}D{\Rightarrow}\tupl{\frak{a,a}}{\in}D$.
Suppose $\tupl{\frak{a,b}}{\in}D$. Since $D$ is symmetric we have $\tupl{\frak{b,a}}{\in}D$.
Since $D$ is not contained in $\tupl{\frak{b,b}}^\times$ there exists $\tupl{\frak{r,s}}{\in}D$ such that
$\frak{b{\notin}\{r,s\}}$. We have\hfil\penalty-10000
$\frak{\tupl{a,a}{=}\tupl{a,b}{\cdot}\tupl{r,s}{\cdot}\tupl{b,a}}{\in}D$.\eod
\vskip5pt
Remark. A similar statement is true without the assumption of symmetry.
Only the case $(\s c)$ should be modified: $D{=}A{\times}B$ for $A,B{\subset}T$.
This can be used in the description of the $\times$-reminder of a subsemigroup of a $\times$-acting group.
We do not need this in this article.
\vskip5pt
We use the following simple observation of \cite[subsection 10]{Ge09}.
A set $D$ of ordered pairs is called 2\it-narrow\rm,
if the following equivalent conditions hold:\hfil\penalty-10000
--- $D$ contain no 3-matching, i.e, a triple whose both projections are triples;\hfil\penalty-10000
--- $D$ is contained in the union of two fibers.

A $\times$-action $\mu:G{\times}T\to T$ of a locally compact group $G$ is called \it elementary \rm if $\partial_\mu G$
is 2-narrow.

\begin{prop}\label{nonElem}
If a $\times$-action $\mu:G{\times}T\to T$ of a locally compact group is not elementary then
each point of $\bold\Lambda G$ is non-isolated in $\bold\Lambda G$ and $\partial G{=}(\bold\Lambda G)^2$.
\end{prop}
\it Proof \rm follows immediately from \ref{xRemainder},
\ref{limSet}, and \ref{nonIsol}.\eod
\vskip5pt
\bf Corollary\sl. A $\times$-action $G{\curvearrowright}T$ possesses fixed points
if and only if $\partial G$ has type $(\s{a})$ or $(\s{b})$ of \ref{xRemainder}.
\vskip3pt
\it Proof\rm.
Let $\partial G{=}\{\frak{\tupl{p,q},\tupl{q,p}}\}$ with $\frak{p{\ne}q}$.
We will prove that the action fixes $\frak p$ and $\frak q$.
Indeed, the set $\{\frak{p,q}\}$ is invariant. Suppose that some $g{\in}G$ transposes $\frak p$ and $\frak q$.
If $h{\in}G$ is sufficiently close to the cross $\tupl{\frak{p,q}}^\times$ then $h^{-1}gh$ is
close to the cross $\tupl{\frak{p,p}}^\times$ which is impossible.

If $\partial G$ has type $(\s b)$ then $\frak p$ is the unique fixed point.

On the other hand if $\partial G$ contains $L^2$ with $|L|{\geqslant}2$
and $\frak p{\in}T$ then
there is a limit cross $\tupl{\frak{r,a}}^\times$ with $\frak{r{\ne}p}$ and $\frak{a{\ne}p}$.
A homeomorphism close to such cross can not fix the point $\frak p$.\eod
\vskip3pt
A $\times$-action with a unique limit cross is called \it parabolic\rm.
The limit cross of a parabolic action has the form $\tupl{\frak{p,q}}^\times$ (so it is of type $(\s b)$)
and the point $\frak p$ is fixed.

Every noncompact locally compact group $G$ possesses a parabolic action:
it is the action on the one-point compactification.
% [cones] @tag_4AB7D4F4:
\subsection{Cones and perspectivity}\label{cones}
Let $\mu:L{\times}X\to Y$ be a $\times$-action.
\deff{ \ref{cones}. An unbounded set $S{\subset}L$ is a \it cone \rm
if
$\partial_1S\cap\overline{\mu(S{\times}\partial_0S)}=\varnothing$.}
An unbounded subset of a cone is a cone. Since continuous maps take
closure to closure, the following immediately follows:
\begin{prop}\label{conPreim} Let $\mu_\iota:L_\iota{\times}X_\iota\to Y_\iota$ be $\times$-actions for $\iota\in\{0,1\}$
and let $\tupl{\alpha,\beta,\gamma}$ be a morphism $\mu_0\to\mu_1$.
If $S{\subset}L_0$ and $\alpha S$ is a cone then $S$ is a cone.\eod\end{prop}

\bf Lemma\sl. If $S$ is a cone for a $\times$-action $\mu:L{\times}X\to Y$ then $|\partial_0S|{=}1$.\hfil\penalty-10000
\sl Proof\rm. Otherwise there are $\frak{\tupl{r,a},\tupl{s,b}}\in\partial S$ with $\frak{r{\ne}s}$.
For arbitrary $b{\in}\s{Loc}_Y\frak b$ there exists $g{\in}S$ such that
$\mu_*g{\subset}\frak s'{\times}Y{\cup}X{\times}b$.
Thus $\mu(g,\frak r){\in}b$.
This implies that $\frak b\in\partial_1S\cap\overline{\mu(S{\times}\partial_0S)}$.\eod
\vskip5pt
The unique point of $\partial_0S$ of a cone $S$ is its \it vertex\rm.
% [conical points] @tag_4DAA1767:
A point $\frak p$ is \it conical \rm if it is a vertex of a cone.
Remark. This notion generalizes the notion of conical point
introduced in \cite{Tu98} and \cite{Bo99} for a convergence action
of a discrete group on a metrisable compactum.

From \ref{conPreim} there follows
\begin{prop}\label{cUnramif} Let $\mu_\iota:L_\iota{\times}X_\iota\to Y_\iota$ be $\times$-actions for $\iota\in\{0,1\}$
and let $\tupl{\alpha,\beta,\gamma}$ be a morphism $\mu_0\to\mu_1$. If $\frak p{\in}X_1$ is a conical point for
$\mu_1$ then the set $\beta^{-1}\frak p$ is a single point which is conical for $\mu_0$.\end{prop}

In \cite[subsection 7]{Ge09} we used another definition of a conical point.
We now prove the equivalence of the two definitions.

\begin{prop}\label{boundedOn2} Let $\mu:L{\times}X\to Y$ be a $\times$-action and let $\beta{=}\{\goth{p,q}\}{\in}\bold\Theta^2X$.
The following properties of a set $S{\subset}L$ are equivalent:\hfil\penalty-10000
$\s a:$ the set $\{S\}\beta{\leftrightharpoons}\{g\beta:g{\in}S\}$ is bounded in $\bold\Theta^2Y$;\hfil\penalty-10000
$\s b:$ $S$ is a union of a bounded set and finitely many cones with vertices in $\beta$.\end{prop}
\it Proof\rm. $\s a{\Rightarrow}b:$
If $\tupl{\frak{r,a}}{\in}\partial S$ then $\frak r{\in}\beta$ since otherwise
$\{S\}\beta$ contains arbitrarily small pairs.
For disjoint closed neighborhoods $P,Q$ of the sets $(\goth p{\times}Y)^{\{\times\}}$
and $(\goth q{\times}Y)^{\{\times\}}$ respectively, $S$ is the union of its intersections
with $\mu_*^{-1}P$ and $\mu_*^{-1}Q$ and a bounded set.

So we can assume that $\partial_0S{=}\{\goth p\}$.

Since $\{S\}\beta$ is bounded in $\bold\Theta^2Y$
it is contained
is a union of finitely many closed subproducts $p{\times}q{\subset}\bold\Theta^2Y$.
So we can assume that $\{S\}(\goth{p,q})\subset p{\times}q$ for disjoint $p,q{\in}\s{Closed}Y$.
Let $\goth a{\in}\partial_1S$. Since $S\goth q$ intersects arbitrary neighborhood of $\goth a$
we have $\goth a{\in}q$. Thus $S$ is a cone.

$\s b{\Rightarrow}a:$ It suffices to consider the case when $S$ is a cone with vertex $\frak p$.
By definition the closed sets $\overline{\mu(S{\times}\frak p)}$ and $\partial_1S$
are disjoint. Let $p,q$ be disjoint closed neighborhoods of these sets.
Then the set $\{s{\in}S:s(\goth{p,q})\notin p{\times}q\}$
 has no limit crosses and therefore is bounded.\eod
\vskip5pt
\bf Corollary\sl. A point $\frak p{\in}X$ is conical for $\mu$ if and only if there exists an unbounded set $S{\subset}L$
such that $\{S\}\{\frak{p,q}\}$ is bounded in $\bold\Theta^2Y$ for every $\frak q{\in}X{\setminus}\frak p$\rm.\eod

So this definition of a conical point is equivalent to that of \cite{Ge09}.
\vskip3pt
% [non-conical pair] @tag_4D895ABA:
Denote by $\s{NC}_\mu$ the set of all non-conical points for a $\times$-action $\mu:L{\times}X\to Y$.

\begin{prop}\label{nonCon2}
For an $\times$-action $\mu:L{\times}X\to Y$,\hfil\penalty-10000
each pair $\beta{\in}\s{NC}_\mu{\leftrightharpoons}\{$non-conical points for $\mu\}$
is \it perspective \rm\sl i.e, the orbit map\hfil\penalty-10000
$L{\ni}g\mapsto g\beta{\in}\b\Theta^2Y$ is proper.\end{prop}
\it Proof\rm. Follows immediately from \ref{boundedOn2}.\eod

\begin{prop}\label{boundedIsPerspective}
For an $\times$-action $\mu:L{\times}X\to Y$,\hfil\penalty-10000
each closed set $B\subset L{\setminus}\partial_1L$
is \it perspective \rm\sl i.e,\hfil\penalty-10000
for every $\b u{\in}\s{Ent}Y$ the set $\s{Big}(B,\b u)$ (see \ref{big}) is bounded.\end{prop}\sl Proof\rm.
by \ref{big} $\partial_0\s{Big}(B,\b u){=}\varnothing$. So it is bounded.\qed
% [image and preimage] @tag_4AB7CBE0:
\subsection{Image and preimage of a $\times$-action}
We say that a map $f:S\to T$ is \it ramified \rm over a point $\frak p{\in}T$
if $|f^{-1}\frak p|{\geqslant}2$.
\begin{prop}\label{xMorph}
Let a locally compact group $G$ act on compactums $X,Y$ and let $f:X\to Y$ be a continuous $G$-equivariant map.
Then\hfil\penalty-10000
$(\s a)$ if $G{\curvearrowright}X$ is 3-proper, $f$ is surjective, and $|Y|{\geqslant}3$ then
$G{\curvearrowright}Y$ is 3-proper;\hfil\penalty-10000
$(\s b)$ if  $G{\curvearrowright}Y$ is 3-proper, $f$ is non-ramified over the limit points, and $|X|{\geqslant}3$ then
$G{\curvearrowright}X$ is 3-proper.
\end{prop}
\it Proof\rm. ($\s a$) follows immediately from the description of $\times$-actions
as Dynkin actions (see \ref{charProp}).

To prove $(\s b)$ consider a limit quasihomeomorphism $s$ for the action $G{\curvearrowright}X$ that is not a homeomorphism.
Since the map $\s{Closed}(X^2)\to\s{Closed}(Y^2)$ induced by $f$ is continuous (see \ref{Vietoris})
and the action $G{\curvearrowright}Y$ is 3-proper,
it maps $s$ to a limit cross $\tupl{\frak{r,a}}^\times$ for $G{\curvearrowright}Y$.
We will prove that the cross $c\leftrightharpoons f^{-1}\frak r{\times}X\cup X{\times}f^{-1}\frak a$ is
a limit quasihomeomorphism for $G{\curvearrowright}X$.
Since $f$ is not ramified over the points $\frak{r,a}$, we have $s{\subset}\tupl{f^{-1}\frak r,f^{-1}\frak a}^\times$.
If follows from \ref{nonIsol} that $s{=}\tupl{f^{-1}\frak r,f^{-1}\frak a}^\times$.\eod
%\subsection{On the uniqueness of morphisms of $\times$-actions}
% [x-property] @tag_4AA83185) end of range
%%%%%%%%%%%%
% [attractorSum] @tag_4DF3262D(
\section{Attractor sum}\label{attSumSect}
% [glueing topology] @tag_4A80F3DF:
\subsection{Gluing topology}
In this subsection we fix a $\times$-action $\mu:L{\times}X\to Y$. For $(g,\goth p){\in}L{\times}X$, $G{\times}P{\subset}L{\times}X$
we put $g\goth p{\leftrightharpoons}\mu(g,\g p)$, $GP{\leftrightharpoons}\mu(G{\times}P)$.
For $K{\times}S{\subset}X{\times}Y$ we denote
\setcounter{equation}0\begin{equation}\label{SKinv}
SK^{-1}\leftrightharpoons\{g{\in}L:S{\cap}gK{\ne}\varnothing\}=((\mu_*)^{-1}((K{\times}S)'))'.\end{equation}
We are going to express the topology of $Y$ in terms of some its restrictions.
\setcounter{prop}1%
\begin{prop}\label{glue}
Let $K{\times}\Lambda\in\s{Closed}(X{\times}Y)$.
Then\hfil\penalty-10000
$\s a:$
If $K{\cap}\partial_0L{=}\varnothing$ then
$\partial_1(FK^{-1}){\subset}F$ for every $F{\in}\s{Closed}Y$;\hfil\penalty-10000
$\s b:$
If $\Lambda'{=}LK{=}LU$ for $U{\in}\s{Loc}_XK$
then\hfil\penalty-10000
$\s{Closed}Y\supset\Cal S\leftrightharpoons\{S:S{\cap}\Lambda{\in}\s{Closed}\Lambda,
S{\cap}\Lambda'{\in}\s{Closed}(\Lambda'),\partial_1(SK^{-1}){\subset}S\}$.\end{prop}\sl Proof\rm.
$(\s a)$.
It follows from \ref{SKinv} that $FK^{-1}$ is a preimage of an $\s{OG}$-closed set.
Hence any limit cross $(\g{r,a})^\times$ for $FK^{-1}$ meets $K{\times}F$. Since $\g r{\notin}K$
we have $\g a{\in}F$.

$(\s b)$ For $\goth p{\notin}S{\in}\Cal S$ we will find $p{\in}\s{Loc}_Y\goth p$ disjoint from $S$.
If $\goth p{\in}\Lambda'$ then $p=\Lambda'{\setminus}S{\cap}\Lambda'$.
Assume $\goth p{\in}\Lambda$.
Since $\g p{\notin}\partial_1(SK^{-1})$ the set $SK^{-1}$
is equicontinuous at $\goth p$ by \ref{eqCont}. Hence, for some $p{\in}\s{Loc}_Y\goth p$,
each $g^{-1}p$, for $g{\in}SK^{-1}$ is small with respect to $\b u\leftrightharpoons\s S^2X{\setminus}(K{\times}U')\in\s{Ent}X$.
Since $\g p{\in}S'{\cap}\Lambda{\in}\s{Open}\Lambda$ we may assume that $p{\cap}\Lambda{\cap}S{=}\varnothing$.
Claim: $p{\subset}S'$.
Indeed if $\goth q{\in}p{\cap}S$ then $\goth q{\notin}\Lambda$ hence $\goth q{=}g\goth r$
for $(g,\g r){\in}L{\times}K$. So $g{\in}SK^{-1}$, $g^{-1}\goth p{\in}U$, $\goth p{\in}\Lambda'$
contradicting with the assumption `$\goth p{\in}\Lambda$'.\eod

\begin{prop}\label{glueUnique}
Let $\Xi{\times}\Lambda$ be a closed invariant subproduct of $X{\times}Y$
with $|\Lambda|{\geqslant}3$ such that $\mu|_{\Xi'{\times}\Lambda'}$ is cocompact.
Then the topology of $Y$ is determined by the topologies of
$\Lambda$ and $\Lambda'$ and by $\nu{\leftrightharpoons}\mu|_{\Xi{\times}\Lambda}$ and $\omega{\leftrightharpoons}\mu|_{\Xi'{\times}\Lambda'}$.\end{prop}
\sl Proof\rm. By \ref{limSet} $\partial_\mu L{=}\partial_\nu L\subset\Xi{\times}\Lambda$.
Let $K$ be a generating compact for $\omega$, i.e, $LK{=}\Lambda'$.
It satisfies the conditions $\s{a,b}$ of \ref{glue}. Hence $\s{Closed}Y{=}\Cal S$.
The definition of $\Cal S$ is a desired expression of the topology of $Y$.\eod
% [particular case] @tag_4DF337B2:
\subsection{Auxiliary sum of $L$ and $\Lambda$}\label{sumWithL}
We transform the definition of $\Cal S$ from \ref{glue} into a definition of a topology.

Let $\nu:L{\times}\Xi\to\Lambda$ be a $\times$-action.
On the direct union $T{\leftrightharpoons}L{\sqcup}\Lambda$ define a topology as follows\hfil\penalty-10000
$\s{Closed}T\leftrightharpoons\{F{\subset}T:F{\cap}L{\in}\s{Closed}L,\ F{\cap}\Lambda{\in}\s{Closed}\Lambda,$ and $\partial_1(F{\cap}L){\subset}F\}$.

Since the operator $\partial_1$ preserves inclusions, the axioms of topology are satisfied.
We denote this construction by $L{+}_\nu\Lambda{\leftrightharpoons}T$.

Since $L$ is locally compact and $\partial($bounded set$){=}\varnothing$ we have $\s{Open}L{\subset}\s{Open}T$.
Thus $L$ is an open subspace of $T$.

Until the end of this subsection we choose and fix $\gamma{\subset}\Lambda$ with $|\gamma|{=}3$.

For a set $s{\subset}\Lambda$
let $\kappa s{\leftrightharpoons}\{g{\in}L:|s{\cap}g\gamma|{>}1\}$,
$\widetilde s{\leftrightharpoons}s{\cup}\kappa s$.
The operators $s\mapsto\kappa s{\subset}L$, $s\mapsto\widetilde s{\subset}T$ preserve inclusion and
commute with the complement $s\mapsto s'$.

If $s{\in}\s{Closed}\Lambda$ then $\widetilde s{\in}\s{Closed}T$.
Indeed we have $\widetilde s{\cap}L{\in}\s{Closed}L$.
If, for $(\goth{r,a}){\in}\partial(\kappa s)$, the attractor $\g a$ were not in $s$
then some $g{\in}\kappa s$ sufficiently close to $(\goth{r,a})^\times$ would map
the set $\gamma{\setminus}\goth r$ of cardinality 2 or 3
to the $\Lambda$-neighborhood $\Lambda{\setminus}s$ of $\goth a$ contradicting with the definition of $\kappa s$.

This proves that $o{\in}\s{Open}\Lambda{\Rightarrow}\widetilde o{\in}\s{Open}T$.
Since $o{=}\Lambda{\cap}\widetilde o$ we have proved the following.\begin{prop}\label{attrSumSubspace}$L$
is an open subspace of $T$ and $\Lambda$ is
a closed subspace of $T$.\hfil\penalty-10000
$\s{Loc}_T\Lambda{=}\{S{\subset}T:S'$ is bounded in $L\}$.\eod\end{prop}
\begin{prop}\label{aSumCompact}$T$ is compact.\end{prop}
\sl Proof\rm.
If an ultrafilter $\Cal F$ on $T$ contains a set bounded in $L$ then it converges to a point in $L$
since $L$ is locally compact. Otherwise it contains $\s{Loc}_T\Lambda$ by \ref{attrSumSubspace}
and converges to a point of $\Lambda$ by \ref{ultra}.\eod
\vskip3pt
For $\b u{\in}\s{Ent}\Lambda$ let $\kappa\b u{\leftrightharpoons}{\cup}\{\s S^2(\kappa s):s{\in}\s{Small}(\b u)\}$,
$\widetilde{\b u}{\leftrightharpoons}{\cup}\{\s S^2\widetilde s:s{\in}\s{Small}(\b u)\}$.
\begin{prop}\label{smallNbhd} Let
$\b w{\in}\Cal T\leftrightharpoons\{\b{v{\cup}\widetilde u:v}{\in}\s{Loc}_{\s S^2L}\b\Delta^2L,\b u{\in}\s{Ent}\Lambda\}$.
Every $\goth p{\in}T$ possesses a $\b w$-small $T$-neighborhood.\end{prop}\sl Proof\rm.
Let $\b{w{=}v{\cup}\widetilde u}$.
If $\goth p{\in}\Lambda$ and $p$ is a $\b u$-small $\Lambda$-neighborhood of $\g p$
then $\widetilde p$ is a $\b{\widetilde  u}$-small $T$-neighborhood of $\goth p$.
If $\g p{\in}L$ and $K$ is a compact $L$-neighborhood of $\goth p$
then any $\b v|_K$-small $K$-neighborhood
of $\goth p$ is also a $\b v$-small $T$-neighborhood of $\goth p$.\eod

We have proved that $\Cal T{\subset}\s{Loc}_{\s S^2T}\b\Delta^2T$.
We are going to proof that $\Cal T$ generates this filter.
\begin{prop}\label{specialSeparation}
Distinct points $\g{p,q}{\in}T$ possess $T$-neighborhoods $p,q$
such that $p{\times}q$ does not meet some
$\b w{\in}\Cal T$.\end{prop}\sl Proof\rm.
Case 1: $\goth{p,q}{\in}L$.
For compact disjoint $L$-neighborhoods $p,q$ of $\goth{p,q}$ respectively let
$\b v\leftrightharpoons\s S^2L{\setminus}p{\times}q\in\s{Loc}_{\s S^2L}\b\Delta^2L$.
There exists $\b u{\in}\s{Ent}\Lambda$ such that
$g\gamma$ does not contain $\b u$-small proper pairs for every $g{\in}p{\cup}q$.
Thus $p{\cup}q$ is disjoint from
$\kappa s$ for any $\b u$-small $s{\subset}\Lambda$.
So $p{\times}q\cap(\b{v{\cup}\widetilde u})=\varnothing$.

Case 2: $\goth p{\in}L$, $\goth q{\in}\Lambda$.
Let $p$ be a compact $L$-neighborhood of $\goth p$ and let $\b u$ be
a $\Lambda$-entourage such that
$g\gamma$ does not contain $\b u$-small proper pairs for every $g{\in}p$.
Let $q$ be a $\b u$-small $\Lambda$-neighborhood of $\goth q$.
The $L$-closed set $\kappa q$ does not meet $p$.
Hence
$\b v\leftrightharpoons\s S^2L{\setminus}p{\times}\kappa q\in\s{Loc}_{\s S^2L}\b\Delta^2L$.
We claim that $p{\times}\widetilde q$ is disjoint from $\b{w{\leftrightharpoons}v{\cup}\widetilde u}$.
If $(g,\goth r){\in}p{\times}q$ then $\{g,\goth r\}$ does not belong to $\widetilde{\b u}$
since $g$ does not belong to $\kappa r$ for any $\b u$-small $r{\subset}\Lambda$.
By the same reason if $(g,h){\in}p{\times}\kappa q$ then the pair $\{g,h\}$ is not in $\kappa\b u$.
It is not in $\b v$ by the choice of $\b v$.

Case 3: $\goth{p,q}{\in}\Lambda$.
Let $\b u$ be a $\Lambda$-entourage with $\{\goth{p,q}\}{\notin}\b u^3$ and let
$p,q$ be closed $\b u$-small $\Lambda$-neighborhoods of $\goth{p,q}$ respectively.
Since the action $\nu$ is continuous the set
$c{\leftrightharpoons}\{(g,h){\in}L^2:\exists\goth r{\in}\gamma:(g\goth r,h\goth r){\in}p{\times}q{\cup}q{\times}p\}$
is $L^2$-closed. Let $\b v{\leftrightharpoons}\s S^2L{\setminus}c\in\s{Loc}_{\s S^2L}\b\Delta^2L$.
Claim: $\widetilde p{\times}\widetilde q$ is disjoint from $\b{w{\leftrightharpoons}v{\cup}\widetilde u}$.

Subcase 1. The set $p{\times}q$ does not contain $\b u$-small pairs by the choice of $\b u$.

Subcase 2. Let $g{\in}\kappa p$ and $\goth r{\in}q$. If $\{g,\goth r\}$ belongs to $\widetilde r$ for
a $\b u$-small $r{\subset}\Lambda$ then $g{\in}\kappa p{\cap}\kappa r{\ne}\varnothing$
hence $p{\cap}r{\ne}\varnothing$.
On the other hand, $\goth r{\in}r{\cap}q{\ne}\varnothing$ so $\{\goth{p,q}\}{\in}\b u^3$.

Subcase 3. Let $g{\in}\kappa p$, $h{\in}\kappa q$.
Since $\gamma{\cap}g^{-1}p{\cap}h^{-1}q{\ne}\varnothing$
the pair $(g,h)$ belongs to $c$ so $\{g,h\}{\notin}\b v$.
If $\{g,h\}{\subset}\kappa r$ for $r{\in}\s{Small}(\b u)$
then, as above, $r$ meets both $p$ and $q$ so $\{\goth{p,q}\}{\in}\b u^3$.\eod

\begin{prop}\label{sumEnt}$T$ is Hausdorff.
The filter $\s{Loc}_{\s S^2T}\b\Delta^2T{=}\s{Ent}T$ is generated by $\Cal T$.\end{prop}\sl Proof\rm.
The first follows from \ref{specialSeparation}.
The second follows from \ref{specialSeparation} and \ref{generatingLoc} applied to $F{\leftrightharpoons}\b\Delta^2T$.\eod
\vskip3pt
It also follows from \ref{sumEnt} that the uniformity $\Cal U$ induced in $L$ by $\s{Ent}T$ is generated
by the set
$\Cal S{\leftrightharpoons}\{\b{v{\cup}\kappa u:v}{\in}\s{Loc}_{\s S^2G}\b\Delta^2G,\b u{\in}\s{Ent}\Lambda\}$.
\begin{prop}\label{sumUnlink} If $\b{u,v}{\in}\s{Ent}\Lambda$ and $\b{u{\bowtie}v}$ then
$\b{\widetilde u{\bowtie}\widetilde v}$ on $T$.\end{prop}\sl Proof\rm.
If $\Lambda{=}a{\cup}b$  where $a{\in}\s{Small}(\b u)$, $b{\in}\s{Small}(\b v)$ then, clearly
$L{=}\kappa a{\cup}\kappa b$ and $T{=}\widetilde a{\cup}\widetilde b$.\eod
\vskip3pt
Suppose now that $X{=}Y$, $L{=}G$ is a group, $\nu$ is a group action.
We consider the action $G{\on}G$ by multiplication from the left.
Let $\mu$ denote the resulting action $G{\on}T$.
The operator $s\mapsto\kappa s$ commutes with the operators of the action $\mu$.
So $\Cal T$ and hence $\s{Ent}T$ is $G$-invariant and $\mu$ acts by homeomorphisms.
Proposition \ref{sumUnlink} implies the following.

\begin{prop}\label{sumWithG}
The action $\mu:G{\on}(G{+}_\nu\Lambda)$ has Dynkin property and hence is a $\times$-action.\eod\end{prop}
We have proved the attractor sum theorem
in the particular case `$\Omega{=}G$'.
It is enough for establishing the Floyd map from any Cayley graph with respect to a finite generating set.
%
% [the theorem] @tag_4DF33779:
\subsection{Existence of attractor sums}\label{sumWithOmega}
We will now prove the attractor sum theorem in full generality.
\begin{prop}\label{attractorSum}
Let
a locally compact group $G$ act on a compactum $\Lambda$ properly on triples
and on a locally compact Hausdorff space $\Omega$ properly and cocompactly.
Then on the disjoint
union $\Omega{\sqcup}\Lambda$
there is a unique compact Hausdorff topology $\tau$
extending the original topologies of
$\Lambda$ and $\Omega$
such that the $G$-action on the space $X{\leftrightharpoons}(\Omega{\sqcup}\Lambda,\tau)$
is proper on triples.\end{prop}

The uniqueness of $\tau$ follows from \ref{glueUnique}, so
it suffices to indicate $\tau$ with the desired properties.

We call the space $X$ the \it attractor sum \rm of $\Lambda$ and $\Omega$. We denote it by
`$\Lambda{+}_G\Omega$' or by $\Omega{+}_G\Lambda$.

\sl Proof of \rm\ref{attractorSum}.
We initially prove the statement for the spaces $\Omega$ of a special type.
We fix a $\times$-action $\nu:G{\times}\Lambda\to\Lambda$ and denote by $T$
the $G$-space $G{+}_G\Lambda{=}G{+}_\nu\Lambda$.
By \ref{aSumCompact} and \ref{sumEnt} $T$ is a compactum and by
\ref{sumWithG} the action $G{\curvearrowright}T$ is 3-proper.

Let $K$ be a compactum, and let first $\Omega$ be $G{\times}K$ with the $G$-action induced
by the trivial action on $K$:
$g(h,\goth p){\leftrightharpoons}(gh,\goth p)$.
This action $G{\on}\Omega$ is proper and cocompact.
We construct the attractor sum $\Omega{+}_G\Lambda$ as follows.

Let $\theta$ be the equivalence relation on the compactum $T{\times}K$ whose classes
are the single points of $\Omega$ and the fibers $\goth p{\times}K$, $\goth p{\in}\Lambda$.
Since $\theta$ is closed the space $X{\leftrightharpoons}(T{\times}K)/\theta$ is a compactum by \ref{Aleksandrov}.
The natural copies of $\Omega$ and $\Lambda$ in $X$ are homeomorphic and $G$-isomorphic
to the corresponding spaces.
We can identify $X$ with $(G{\times}K)+_G\Lambda$.
Indeed the natural projection $\pi:X\to T$ is continuous $G$-equivariant.
By \ref{xMorph}($\s b$) the action
$G{\on}X$ is 3-proper.

Now consider the general case of a proper cocompact action $\omega:G{\times}\Omega\to\Omega$ on a locally compact space.
Let $K$ be a compactum in $\Omega$ such that $\Omega{=}GK$ and let $\theta$
be the kernel of the surjective proper map $\omega|_{G{\times}K}:G{\times}K\to\Omega$.

We claim that the equivalence $\Theta{\leftrightharpoons}\theta{\cup}\b\Delta^2\Lambda$ is
closed in $X^2$.

Let $\theta_1$ be the image of $\theta$ under the map $\pi^2:X^2\to T^2$.
So $\theta_1{=}{\cup}\{gB^2:g{\in}G\}$ where
$B{\leftrightharpoons}\{g{\in}G:K{\cap}gK{\ne}\varnothing\}$.
Since $\pi^2|_{(G\times K)^2}$ is a proper map, $\theta_1$ is closed in $G^2$.
The action $G{\curvearrowright}\Omega$ is proper so
$B$ is bounded and hence perspective (see \ref{boundedIsPerspective}).
So
the set $\overline{\theta_1}{\setminus}\theta_1$, where bar means the closure in $T^2$,
does not contain non-diagonal pairs. Hence $\overline{\theta_1}{=}\theta_1{\cup}\b\Delta^2\Lambda$.
The $\pi^2$-preimage of $\overline{\theta_1}$ is contained
in $(G{\times}K)^2{\cup}\b\Delta^2\Lambda$.
Since it is closed and contains $\theta$ we have $\overline\theta{\subset}(G{\times}K)^2{\cup}\b\Delta^2\Lambda$.
Since $\theta$ is closed in $(G{\times}K)^2$ we have $\overline\theta{\subset}\theta{\cup}\b\Delta^2\Lambda$.

So $\Theta$ is closed.
The space $\Omega{+}_G\Lambda{\leftrightharpoons}((G{\times}K)+_G\Lambda)/\Theta$
is a compactum by \ref{Aleksandrov}.
The action of $G$ on
$(G{\times}K)+_G\Lambda$ is proper on triples, so by
\ref{xMorph}(a) the induced action $G{\on}(\Omega{+}_G\Lambda)$
is proper on triples too.
Since $\omega|_{G{\times}K}$ is proper the topology of $\Omega$ coincides with
the quotient topology, so $\Omega$ is a subspace of $\Omega{+}_G\Lambda$.\qed
\subsection{Implication $\s{RH_{32}\Rightarrow RH_{pd}}$}
Let a discrete group $G$ act on a compactum $T$ properly on triples and cocompactly on pairs.
Let $\widetilde T{\leftrightharpoons}G{+}_GT$.
By \ref{attractorSum} the action $G{\on}\widetilde T$ is proper on triples.
If $W$ is a compact fundamental domain for the action $G{\on}\b\Theta^2T$ then
$\widetilde W{\leftrightharpoons}W{\cup}(1{\times}(\widetilde T{\setminus}1))$ is
a compact fundamental domain for the action
$G{\on}\b\Theta^2\widetilde T$.
Let $\widetilde{\b u}$ be the complement of a compact neighborhood $\widetilde W_1$ of $\widetilde W$
in $\b\Theta^2\widetilde T$.
By \ref{generatingLoc} the $G$-filter generated by $\{\widetilde{\b u}\}$
coincides with $\s{Ent}\widetilde T$.
Let $\widetilde{\b v}$ be an open entourage of $\widetilde T$ such that
$\widetilde{\b v}^2{\subset}\widetilde{\b u}$.
Then there exists a finite set $F{\subset}G$ such that $\widetilde{\b v}'\subset{\cup}F\{\widetilde W_1\}$.
Then $\widetilde{\b v}\supset{\cap}F\{\widetilde{\b u}\}$,
and so $({\cap}F\{\widetilde{\b u}\})^2\subset\widetilde{\b u}$.
So $\widetilde{\b u}$ is a divider (see \ref{expansiveMetric}).

Let $M$ be the set of all non-conical points of $G{\on}\widetilde T$.
It contains $G$ and is $G$-finite by \cite[Section 7.2]{Ge09}.
The restriction $\b u{\leftrightharpoons}\widetilde{\b u}{\cap}\b\Theta^2M$ is
a perspective divider on $M$
by \ref{nonCon2}. So we have the following.
\begin{prop}\label{relHyperbolicityOf32}Let $G{\curvearrowright}T$ be
an $\s{RH_{32}}$-action and the $G$-set $M$ of nonconical points of $\widetilde T$ is connected
then the uniformity of $\widetilde T$ induces a relative hyperbolicity on $M$.\end{prop}

If $G$ is finitely generated then $G{\on}G$ is connected since there exists a locally finite Cayley graph.
In \cite{GP10} we will prove a theorem which implies that $M$ is always connected.
So in all cases we have a connected $G$-set
$M$ which admits a relative hyperbolicity.

The implication $\s{RH_{32}{\Rightarrow}RH_{pd}}$ is proved.
% [attractorSum] @tag_4DF3262D)
\section{Remarks}
\subsection{Bowditch completion coincides with the completion
with respect to the relatively hyperbolic uniformity}\label{BoCmpl}
Starting from an $\s{RH_{fh}}$-action $G{\curvearrowright}\Gamma$
of a finitely generated group $G$, B. Bowditch \cite{Bo97} constructed
a compactum $B$ and an action $G{\curvearrowright}B$
such that $B{\supset}\Gamma^0$ (equivariantly),
the spaces of the form $(B{\setminus}v)/\s{St}_Gv$
($v{\in}\Gamma^0$) are compact,
the action $G{\curvearrowright}B$ has the convergence property,
and each point of $B{\setminus}\Gamma^0$ is a conical limit point.
C. Hruska \cite{Hr10} extended this result to countable groups $G$.
An intermediate step is a construction of a ``proper hyperbolic
length space'' $X$ with a ``geometrically finite'' action $G{\curvearrowright}X$ by isometries.
This $X$ is not canonical. However Bowditch proved
\cite[Theorem 9.4]{Bo97} that the space $B$ with the indicated properties is unique.

We are going to explain that Bowditch's space $B$
 can be obtained directly as a completion
with respect to the visibility
described in the beginning of Section \ref{altRelHyp}.
Since we do not need an intermediate metric space $X$
we need no restriction on the cardinality of $G$.

Let $\Gamma$ be as in the definition $\s{RH_{fh}}$.
By \ref{fineHypIsAltHyp} $\Gamma$ is alt-hyperbolic so
the uniformity $\Cal U$ on $M{\leftrightharpoons}\Gamma^0$, derermined by the divider $\b u_E$ of
\ref{altRelHyp} is a visibility.
Clearly, it is exact (see \ref{Samuel})
thus the completion map $\iota_{\Cal U}:M\to B{\leftrightharpoons}\overline M$
is injective. We identify $\Gamma^0$ with $\iota_{\Cal U}\Gamma^0$.
By \ref{pd2fh} the action $G{\curvearrowright}B$ has property $\s{RH_{32}}$.
By \cite[Main Theorem]{Ge09} every non-conical limit point of $B$ is bounded parabolic.

Claim: a vertex $v{\in}\Gamma^0$ can not be conical.

Proof:
let $w$ be a vertex in $\Gamma^0{\setminus}v$
joined with $v$ by an edge
and let $S$ be an infinite subset of $G$.
The set $S\{v,w\}$ is unbounded in $\b\Theta^2B$ by perspectivity and we apply
the corollary of \ref{boundedOn2}.\qed

So each vertex is either isolated or bounded parabolic on $B$ depending on
the finiteness of its $G$-stabilizer.

Finally we prove that every $\goth p{\in}B{\setminus}\Gamma^0$ is conical.
Recal that $\goth p$ is a minimal $\Cal U$-Cauchy filter (see \ref{Samuel}) different from
$\s{Loc}_{\Cal U}v$ for each $v{\in}\Gamma^0$.
Since $B$ is Hausdorff, by the generalized Karlsson lemma \ref{KarlssonLemma},
for every finite set $V{\subset}\Gamma^0$,
there exists a finite $E(V){\subset}\Gamma^1$
that ``overshadows'' $\goth p$ from $V$:
for each sufficiently small $A{\in}\goth p$ every geodesic
segment joining a vertex in $V$ with a vertex in $A$ pass through an edge in $E(V)$.
In particular, sufficiently small sets in $\goth p$ are disjoint from $V$.
This allows us, by an easy induction, to construct an infinite geodesic ray
starting from a given vertex $v{\in}\Gamma^0$ and converging to $\goth p$.
Since $\Gamma^1/G$ is finite there exists an orbit that has
in the ray infinitely many edges. Let $\{g_ie:i=1,2,\dots\}$ be these edges.
The set $\{g_i^{-1}\{v,\goth p\}:i=1,2,\dots\}$ is bounded in $\b\Theta^2B$
since the pairs in this set do not belong to $\overline{\b u}$ where $\b u{=}\root3\of{\b u_e}$.
So $\goth p$ is conical (see the corollary of \ref{boundedOn2}).
\vskip3pt
Thus $B$ is the Bowditch completion of $\Gamma$.
To prove Bowditch's uniqueness theorem 9.4 without restriction on cardinality one has
to verify that different connecting structures (\ref{iDivider}) determine
the same uniformity. The reader can verify that this would follow from a version of the
generalized Karlsson lemma \ref{KarlssonLemma} \bf for quasi-geodesics \rm instead of geodesics.
The proof of this version is actually the same: initially one proves the original Karlsson's
version and then apply the same sequence of estimates.
% [geometric] @tag_4E03ABB1:
\subsection{Geometric actions}\label{geometric}
Let a discrete finitely generated group $G$ act on a compactum $T$ properly on triples.
Let $\Cal U$ be the uniformity on $G$ induced by the uniformity of $G{+}_GT$.

We say that the action $G{\on}T$ is \it geometric \rm with respect to a finite generating set $S$
if $\Cal U$ is a visibility (see \ref{visibility}) on the Cayley graph $\Gamma$ with respect to $S$.

By \ref{visIsDynkin} every geometric action has convergence property.
By Karlsson lemma \cite{Ka03} every action of a f.g. group on its Floyd completion with
respect to any Floyd function on a locally finite Cayley graph is geometric.

It follows directly from the definition that the
quotient of a geometric action is geometric
and that the inverse limit of geometric actions is geometric.

By \ref{mapThrm} the action of a relatively hyperbolic group on its Bowditch boundary is geometric.
In particular the action of a hyperbolic group on its Gromov boundary is geometric.
If $H$ is a hyperbolic subgroup of hyperbolic group $G$ such that the inclusion $H{\subset}G$
induces the ``Cannon-Thurston map'' $\partial_\infty H\to\partial_\infty G$ in
particular when $H$ is normal in $G$ of infinite index then
the (non-geometrically finite) action
$H{\on}\b\Lambda H$ is geometric (as a quotient of the geometric action $H{\on}\partial_\infty H$)
\cite{Mi98}.

\sc Question 1\rm. Does geometricity depend on the choice of finite generating set?

\sc Question 2\rm. Does there exist a non-geometric convergence action of a f.g. group?
\subsection{Convergence actions on totally disconnected spaces}\label{totDisc}
We will show that every convergence action of a discrete f.g. group on a totally disconnected
compactum is geometric. It follows from the following proposition (compare with \cite[Subsectin 7.1,
Proposition E]{Ge09})

\begin{prop}\label{tDisconnect}
Let a finitely generated group $G$ act on a totally disconnected compactum $T$
with the convergence property
and let $\Cal B$ be the Boolean algebra $\s{Open}T{\cap}\s{Closed}T$.
Then the following properties are equivalent:\hfil\penalty-10000
$\s{a}:\b\Theta^2T/G$ is compact;\hfil\penalty-10000
$\s{b}:$ $\Cal B$ is finitely generated as an $G$-boolean algebra.
\end{prop}
\it Proof\rm.
The finite subalgebras of $\Cal B$ correspond to the partitions of $T$ into
finitely many
open subsets. The pieces of a partition are the atoms of the subalgebra.
For a finite subalgebra $B$ of $\Cal B$
the set $\b u_B{\leftrightharpoons}{\cup}\{\s S^2o:o{\in}\s{At}B{\leftrightharpoons}\{$the
atoms of $B\}\}$ is an entourage of $T$.
Moreover, it is a divider (for $F{=}\{1_G\}$, see definition \ref{expansiveMetric}).
Every entourage contains an entourage of the form $\b u_B$.
It is easy to see that $\b u_B$ generates the $G$-filter $\s{Ent}T$ if and only if
$\Cal B{=}\sum G\{B\}$.\qed
\vskip3pt
\bf Corollary\sl. Every $\times$-action of a discrete group $G$ on a totally disconnected
space is an inverse limit of 2-cocompact actions\rm.

\vskip3pt
The following fact is actually proved by J. Stallings \cite{St71}.

\begin{prop}The action of a finitely generated group $G$ on its Freudenthal
completion $\s{Fr}G$ is geometric.
So it has the convergence property.\end{prop}
\it Proof\rm. See \ref{visibility}.\qed
\vskip3pt
It now follows from \ref{tDisconnect} and the discussion of \ref{visibility}
that the action $G{\curvearrowright}\s{Fr}G$ is the inverse limit of 2-cocompact actions.
It is 2-cocompact if and only if it is relatively hyperbolic;
the former is true if and only if the group is accessible in the sense of Dunwoody.
See \cite[Corollary IV.7.6]{DD89} for further details.

Claim: this action $G{\curvearrowright}\s{Fr}G$ is the ``initial'' (see \ref{visibility}) in the class of the convergence actions of $G$ on totally disconnected
spaces.

\it Sketch of the proof\rm.
Let $\s{AI}(G)$ denote the Boolean algebra of subsets of $G$ almost invariant from the right.
That is $\s{AI}(G){\leftrightharpoons}\{A{\subset}G:\forall g{\in}G\ |A{\setminus}Ag|{<}\infty\}$.
It is well-known (see \cite{St71}) that the space $\s{Fr}G$ can be identified
with the space of maximal ideals $\s{Spec\,AI}(G)$ of the algebra $\s{AI}(G)$.
The group $G$ can be identified with an open discrete invariant subspace of $\s{Spec\,AI}(G)$:
to an element $g{\in}G$ there corresponds the ideal $I_g{\leftrightharpoons}\{A{\in}\s{AI}(G):g{\notin}A\}$.

Consider a convergence action $G{\curvearrowright}T$ on a totally disconnected compactum.
The attractor sum space $G{+}T$ (see \ref{sumWithG})
is also totally disconnected, so, by the attractor sum theorem \ref{attractorSum}
(actually we need only the particular case \ref{sumWithG}),
we can assume that $G{\subset}T$.
Moreover we can assume that $G$ is dense in $T$. We say in this case that
$T$ is a $\times$\it-completion \rm of $G$.

Let $\Cal B{\leftrightharpoons}\s{Open}T{\cap}\s{Closed}T$.
We will indicate a $G$-equivariant homomorphism $\varphi:\Cal B\to\s{AI}G$.
By Stone's contravariant equivalence of the category of Boolean algebras and
the category of totally disconnected
Hausdorff compacta it indices the continuous $G$-equivariant map
$\s{Spec\,AI}(G)\to T{\simeq}\s{Spec}\Cal B$.

We put $A^\varphi{\leftrightharpoons}\{g{\in}G:g^{-1}{\in}A\}$.
So it suffices to prove that if $A{\in}\Cal B$
then $A{\cap}G$ is almost invariant \bf from the left\rm.
Suppose that
the set $S\leftrightharpoons A{\cap}G\setminus g(A{\cap}G)$ is infinite
for some $g{\in}G$. Let $(\g{p,q})^\times$ be a limit cross for $S$.
Since $A$ is closed we have $\g q{\in}A$.
On the other hand $(\g{p,q})^\times$ is a limit cross for $g(A'{\cap}G)$.
It follows from \ref{eqCont} that $\partial_1(g(A'{\cap}G)){=}\partial_1(A'{\cap}G)$.
Since $A'$ is closed we have $\g q{\in}A'$. The contradiction proves the claim.\qed

% [bliography] @tag_4DF16318:
\end{document}